\documentclass[11pt]{article}

\usepackage[T1]{fontenc}
\usepackage[a4paper,top=27mm,bottom=28mm,inner=30mm,outer=30mm]{geometry}
\usepackage{amsmath,amssymb,amsthm,mathtools,bm}
\usepackage{enumitem}
\usepackage{needspace}
\usepackage[hidelinks,hyperfootnotes=false]{hyperref}
\usepackage{xcolor}
\newif\ifshowchanges
\showchangesfalse
\definecolor{revisionblue}{RGB}{0,65,175}
\newcommand{\rev}[1]{\ifshowchanges{\color{revisionblue}#1}\else#1\fi}
\newenvironment{revision}
  {\begingroup\ifshowchanges\color{revisionblue}\fi}
  {\endgroup}
\usepackage{microtype}
\usepackage{mathptmx}
\usepackage{fancyhdr}

\hypersetup{
	pdftitle={Critical Optimizers for the Finite-Rank Fractional Lieb-Thirring Inequality},
	pdfsubject={Critical finite-rank optimizers, Birman-Schwinger eigenvalues and profile decompositions}
}

\newtheorem{theorem}{Theorem}[section]
\newtheorem{lemma}[theorem]{Lemma}
\newtheorem{proposition}[theorem]{Proposition}
\newtheorem{corollary}[theorem]{Corollary}

\theoremstyle{definition}

\theoremstyle{remark}

\newtheorem*{unnumberedremark}{Remark}

\newcommand{\R}{\mathbb R}
\newcommand{\C}{\mathbb C}
\newcommand{\N}{\mathbb N}
\newcommand{\one}{\mathbf 1}

\newcommand{\ipHs}[2]{\left\langle #1,#2\right\rangle_{\dot H^s}}
\newcommand{\normHs}[1]{\left\|#1\right\|_{\dot H^s}}
\newcommand{\weak}{\rightharpoonup}
\newcommand{\Tr}{\operatorname{Tr}}

\newcommand{\supp}{\operatorname{supp}}
\newcommand{\diag}{\operatorname{diag}}

\newcommand{\cT}{\mathcal T}

\title{\MakeUppercase{Critical Optimizers for the Finite-Rank\\Fractional
		Lieb--Thirring Inequality}}
\author{Zhuocheng Li\qquad Xiaoyu Zeng\\[0.3em]
\normalsize\itshape Center for Mathematical Sciences, Wuhan University of Technology,\\
\normalsize\itshape Wuhan 430070, P.R. China}

\date{October 24, 2026}

\begin{document}
	\maketitle
	\thispagestyle{plain}
\begingroup
\renewcommand{\thefootnote}{}
\footnotetext{E-mail: \href{mailto:lzc2023203958@whut.edu.cn}{lzc2023203958@whut.edu.cn} (Z.C. Li); \href{mailto:xyzeng@whut.edu.cn}{xyzeng@whut.edu.cn} (X.Y. Zeng).}
\endgroup
	
	\begin{abstract}
\begin{revision}
We study optimal potentials for the finite-rank fractional Lieb--Thirring
problem at the eigenvalue-counting endpoint, for \(0<s<1\) and \(d>2s\).
Using the zero-energy Birman--Schwinger formulation, we prove a profile
decomposition of normalized optimizing sequences with exact splitting of
rank and potential mass, and strong approximation in \(L^{d/(2s)}\).
Strict subadditivity yields existence and compactness modulo translations
and dilations. At rank two, the optimal constant equals the rank-one
constant, the supremum is not attained, and every normalized optimizing
sequence splits into two Sobolev potentials. We also prove strict increase
of the minimal potential masses and a uniform spectral gap for optimizing
sequences. For each optimizer, we derive an Euler--Lagrange system and
establish boundedness, local H\"older regularity and
\(V(x)=|x|^{-4s}(c_\infty+o(1))\) as \(|x|\to\infty\), with
\(c_\infty\geq0\). The conformal orbit of every optimizer contains a
representative with a strictly positive asymptotic coefficient.
\end{revision}
\end{abstract}
	
	
	
	\setcounter{tocdepth}{2}
	\begingroup
	\small
	\tableofcontents
	\endgroup
	\clearpage
	
	\section{Introduction}\label{sec:introduction}

Lieb--Thirring inequalities control moments of the negative spectrum of
Schr\"odinger operators by integral norms of the potential. Their original
applications concern stability of matter~\cite{LT75,LT76}; broader accounts
of the spectral and variational questions can be found in
\cite{FrankSurvey,FLW,LiebSeiringer}.
This paper concerns the eigenvalue-counting endpoint for the fractional
Schr\"odinger operator. Throughout, we assume
\begin{equation}\label{eq:critical-regime}
0<s<1,\qquad \gamma=0,\qquad d>2s.
\end{equation}
Write \(N_-(H)\) for the number of strictly negative eigenvalues of a
self-adjoint operator \(H\), counted with multiplicity.
The fractional Cwikel--Lieb--Rozenblum inequality
\cite[Example~3.3]{FrankCLR} states that
\[
N_-\bigl((-\Delta)^s-V\bigr)
\leq L_{0,d,s}\int_{\R^d}V^{d/(2s)},
\qquad 0\leq V\in L^{d/(2s)}(\R^d).
\]
The local eigenvalue-counting estimate originates in the work of
Cwikel, Lieb and Rozenblum~\cite{Cwi77,Lie76,Roz72}.

Finite-rank optimization connects the rank-one problem with the full
spectral inequality. For the local operator \(-\Delta\), this approach was
developed by Frank, Gontier and Lewin~\cite{FGL21a,FGL}. In particular,
their analysis of the critical case \(\gamma=0\), \(d\geq3\), uses a
variational principle for the eigenvalues of the Birman--Schwinger operator.
They establish a profile decomposition, a compactness criterion, and a
spectral gap and Euler--Lagrange equation for optimizers; see
\cite[Section~3, Theorems~7 and~9, and Corollary~8]{FGL}.
\begin{revision}
The present work extends this critical framework to the fractional
operator \((-\Delta)^s\), with \(0<s<1\) and \(d>2s\).
The profile decomposition, compactness criterion and optimizer system
below are the fractional counterparts of those results. We adapt the
spectral arguments of \cite[Sections~8 and~9]{FGL}, using the fractional
Sobolev profile decomposition and a Riesz-potential argument for regularity
and Kelvin asymptotics. We also establish the rank-two nonattainment
statement corresponding to \cite[Theorem~10]{FGL}; its proof uses the
strict interaction between the positive and negative parts in the
fractional Dirichlet form.
\end{revision}

For fractional operators and positive Riesz exponents
\(\gamma>\max\{0,1-d/(2s)\}\), Chen and Yang~\cite{CY26} prove existence
of optimizers for every finite rank and study their spectral gap,
Euler--Lagrange equation, regularity and decay. Under the additional
condition \(\gamma+d/(2s)\geq3+4s/d\), they prove precompactness of
normalized maximizing sequences modulo translations.
Their positive-moment results do not include the endpoint considered here.

At \(\gamma=0\), strictly negative eigenvalues are counted without an
energy weight. The resulting counting quotient cannot attain its
supremum: whenever the negative-eigenvalue count is positive, a sufficiently
small reduction of the coupling preserves that count while decreasing
the potential mass. The same observation applies to truncated counting.
This motivates a threshold formulation in terms of the zero-energy
Birman--Schwinger operator. Its eigenvalues strictly larger than one count
the negative eigenvalues of the Schr\"odinger operator, whereas the
eigenvalue one corresponds to a nontrivial solution of
\(( -\Delta)^s u=Vu\) in \(\dot H^s(\R^d)\), which need not belong to
\(L^2(\R^d)\). Moreover, the critical scaling
\[
V_t(x):=t^{2s}V(tx),\qquad t>0,
\]
preserves both the negative-eigenvalue count and \(\int V^{d/(2s)}\).
Fixing this integral therefore does not exclude concentration at different
spatial scales. We study the \(N\)-th Birman--Schwinger eigenvalue and
analyze optimizing sequences modulo translations and dilations.
The relation to truncated eigenvalue counting is made explicit below.

At rank one, the problem is equivalent to the sharp fractional Sobolev
inequality, whose optimizers are known from Lieb's work on the sharp
Hardy--Littlewood--Sobolev and related inequalities~\cite{Lieb83}.

We establish a profile decomposition for every normalized critical
optimizing sequence, with exact splitting of the rank and potential mass
and strong approximation in \(L^{d/(2s)}\) by translated and dilated
profiles. Here the rank of a profile counts its Birman--Schwinger
eigenvalues at or above one, including multiplicities.
\begin{revision}
At rank one, the decomposition recovers the compactness of the Sobolev
problem modulo translations and dilations. At rank two, we prove
\(a_2=2a_1\), nonattainment, and decomposition into exactly two rank-one
Sobolev potentials. At ranks \(N\geq3\), existence and compactness follow
under strict subadditivity. The decomposition also implies
\(a_{N+1}>a_N\) for every \(N\), and consequently a uniform gap below
the threshold for normalized optimizing sequences.
\end{revision}
We also prove a spectral gap below the Birman--Schwinger eigenvalue one
and derive the Euler--Lagrange system for optimizers. Every optimizer is
bounded and locally H\"older continuous and satisfies
\[
V(x)=|x|^{-4s}\bigl(c_\infty+o(1)\bigr),\qquad c_\infty\geq0,
\quad |x|\to\infty.
\]
Starting from any given optimizer, a translation followed by Kelvin
inversion produces an optimizer for the same constant with a strictly
positive asymptotic coefficient.

The proof combines a vector-valued homogeneous Sobolev profile
decomposition~\cite{Gerard,PalatucciPisante} with Birman--Schwinger
min--max arguments. The spectral gap excludes excess threshold
multiplicity and is a key ingredient in obtaining the exact rank identity.
Equality in the mass bounds then yields strong approximation of the
potentials. One-sided perturbation theory~\cite{Kato} and convex
separation give the Euler--Lagrange system, while improved integrability
and the Kelvin transform yield regularity and decay.

\subsection{Variational formulation}

We first fix the operator and spectral notation for the critical variational problem.
For \(0<s<1\), the fractional Laplacian is defined on \(\mathcal S(\R^d)\) by
\[
    \widehat{(-\Delta)^sf}(\xi)=|\xi|^{2s}\widehat f(\xi).
\]
We write
\[
    \normHs f^2=\int_{\R^d}|\xi|^{2s}|\widehat f(\xi)|^2\,d\xi,
    \qquad
    \|f\|_{H^s}^2=\|f\|_2^2+\normHs f^2.
\]
When \(d>2s\), the homogeneous space \(\dot H^s(\R^d)\) embeds into
\(L^{2_s^*}(\R^d)\), where \(2_s^*:=2d/(d-2s)\).
The analytic estimates used below are collected in Section~\ref{sec:basic-results}.

For \(V\geq0\), the quadratic form
\(\normHs u^2-\int_{\R^d}V|u|^2\), with domain \(H^s(\R^d)\), defines the self-adjoint operator
\[
    H_{s,V}:=(-\Delta)^s-V.
\]
\begin{revision}
For \(V\in L^{d/(2s)}(\R^d)\), truncation and the critical Sobolev
inequality give, for every \(\varepsilon>0\),
\[
 \int V|u|^2\leq\varepsilon\|u\|_{\dot H^s}^2
       +C_\varepsilon\|u\|_2^2,\qquad u\in H^s(\R^d).
\]
The form is therefore closed and bounded from below; we use the standard
form and min--max framework of \cite{RS75,RS78}.
\end{revision}
We denote the negative eigenvalues, counted with multiplicity, by
\[
    \lambda_1(H_{s,V})\leq\lambda_2(H_{s,V})\leq\cdots\leq0,
\]
with the convention that \(\lambda_j=0\) for \(j>N_-(H_{s,V})\).

Under \eqref{eq:critical-regime}, put
	\begin{equation}\label{eq:critical-exponents}
		q_c:=\frac d{2s},
		\qquad
		q_c':=\frac d{d-2s}.
	\end{equation}
	For \(0\leq V\in L^{q_c}(\R^d)\), define the compact nonnegative operator
	\begin{equation}\label{eq:critical-birman-schwinger}
		K_{s,V}(0)
		:=
		V^{1/2}(-\Delta)^{-s}V^{1/2}
	\end{equation}
	on \(L^2\), or equivalently the compact operator \(\cT_V\) on \(\dot H^s\) determined by
	\begin{equation*}
		\ipHs f{\cT_Vg}
		=
		\int_{\R^d}Vf\overline g.
	\end{equation*}
	The nonzero eigenvalues of the two realizations coincide; write \(\mu_1(V)\geq\mu_2(V)\geq\cdots\geq0\) for their common values, counted with multiplicity.

The critical finite-rank constant is
	\begin{equation}\label{eq:critical-functional}
		\ell_{0,d,s}^{(N)}
		:=
		N\sup_{\substack{0\leq V\in L^{q_c}\\V\neq0}}
		\frac{\mu_N(V)^{q_c}}
		{\displaystyle\int_{\R^d}V^{q_c}}.
	\end{equation}
	Equivalently, define \(a_N\) by the variational characterization
	\begin{equation}\label{eq:critical-minimization}
		a_N
		:=\frac{N}{\ell_{0,d,s}^{(N)}}
		=
		\inf\left\{
		\int_{\R^d}V^{q_c}:
		V\geq0,\ \mu_N(V)\geq1
		\right\}.
	\end{equation}
An \emph{optimizer} for \(\ell_{0,d,s}^{(N)}\) is a nonzero nonnegative
potential \(V\in L^{q_c}(\R^d)\) attaining the supremum in
\eqref{eq:critical-functional}. It is \emph{normalized} if \(\mu_N(V)=1\);
equivalently, a normalized optimizer attains the infimum in
\eqref{eq:critical-minimization}. Unless stated otherwise, the term
optimizer refers to \(\ell_{0,d,s}^{(N)}\), not to the truncated counting
constant introduced below. Here \(N\) specifies the selected eigenvalue
index; it does not assert that the Birman--Schwinger operator has finite
rank or that the potential has exactly \(N\) strictly negative eigenvalues.

Since \(\mu_N(cV)=c\mu_N(V)\) for \(c>0\), amplitude scaling leaves the quotient in \eqref{eq:critical-functional} unchanged.
Thus a maximizing sequence with positive \(\mu_N\) may be normalized by \(\mu_N(V_n)=1\); equivalently, a minimizing sequence for \(a_N\) may be chosen with equality in the constraint \(\mu_N(V)\geq1\).
We refer to \(\int V^{q_c}\) as the mass of the potential; this convention
also applies to the profiles below. This amplitude normalization does not
remove the spatial dilation symmetry.
The action of translations and dilations and its invariances are
	\begin{equation*}
		V_{t,x_0}(x)
		:=
		t^{2s}V\bigl(t(x-x_0)\bigr),
		\qquad
		\mu_j(V_{t,x_0})=\mu_j(V),
		\qquad
		\int V_{t,x_0}^{q_c}=\int V^{q_c}.
	\end{equation*}
We use precompactness modulo translations and dilations in the following
sequential sense: every subsequence admits a further subsequence which,
after suitable translations and dilations, converges strongly in \(L^{q_c}\).

The full fractional CLR constant is
	\begin{equation*}
		L_{0,d,s}
		:=
		\sup_{\substack{0\leq V\in L^{q_c}(\R^d)\\V\neq0}}
		\frac{N_-(H_{s,V})}{\displaystyle\int_{\R^d}V^{q_c}}.
	\end{equation*}
	There is a related finite-rank truncation of the counting problem,
	\begin{equation*}
		\widehat L_{0,d,s}^{(N)}
		:=
		\sup_{\substack{0\leq V\in L^{q_c}(\R^d)\\V\neq0}}
		\frac{\min\{N,N_-(H_{s,V})\}}
		{\displaystyle\int_{\R^d}V^{q_c}},
	\end{equation*}
	and the Birman--Schwinger principle gives
	\begin{equation}\label{eq:critical-truncated-relation}
		\widehat L_{0,d,s}^{(N)}
		=
		\max_{1\leq k\leq N}\ell_{0,d,s}^{(k)}.
	\end{equation}
	Thus \(\ell_{0,d,s}^{(N)}\) concerns the \(N\)-th Birman--Schwinger eigenvalue, whereas \(\widehat L_{0,d,s}^{(N)}\) concerns the negative-eigenvalue count truncated at \(N\). The constants \(\ell_{0,d,s}^{(k)}\) retain the individual rank information needed for the profile decomposition. Equality \(\widehat L_{0,d,s}^{(N)}=\ell_{0,d,s}^{(N)}\) holds precisely when the maximum in \eqref{eq:critical-truncated-relation} is attained at \(k=N\).
The identity \eqref{eq:critical-truncated-relation} is justified after the Birman--Schwinger principle in Section~\ref{sec:basic-results}.

For completeness, the connection with the rank-one Sobolev problem can
be written explicitly. Set
\[
 S_{d,s}:=\inf_{0\neq u\in\dot H^s(\R^d)}
 \frac{\|u\|_{\dot H^s}^2}{\|u\|_{2_s^*}^2}.
\]
The Rayleigh principle and H\"older duality give
\[
 \sup_{0\neq V\geq0}\frac{\mu_1(V)}{\|V\|_{q_c}}
 =\sup_{0\neq u\in\dot H^s}
 \frac{\|u\|_{2_s^*}^2}{\|u\|_{\dot H^s}^2}
 =S_{d,s}^{-1},
 \qquad a_1=S_{d,s}^{q_c},
\]
where the first supremum is over \(V\in L^{q_c}\).
Equality in H\"older's inequality associates a Sobolev optimizer \(u\)
with a potential proportional to \(|u|^{4s/(d-2s)}\), which can be
normalized by \(\mu_1(V)=1\).

\begin{revision}
The equality cases in the sharp Sobolev inequality~\cite{Lieb83} show
that every normalized rank-one optimizer is a translate and a dilation of
\[
 Q_{d,s}(x):=A_{d,s}(1+|x|^2)^{-2s},
 \qquad A_{d,s}>0\text{ chosen so that }\mu_1(Q_{d,s})=1.
\]
Indeed, for a normalized rank-one optimizer \(V\), take a
\(\dot H^s\)-normalized eigenfunction \(u\) at \(\mu_1(V)=1\). Then
\[
 1=\int V|u|^2\leq S_{d,s}\|u\|_{2_s^*}^2\leq1,
\]
so equality holds both in H\"older and in the sharp Sobolev inequality.
\end{revision}

	\subsection{Main results}\label{sec:main-critical}

We now assume \eqref{eq:critical-regime}.
The quantities \(a_N\) defined in \eqref{eq:critical-minimization} are nondecreasing and subadditive; these elementary properties, established in Lemma~\ref{lem:critical-constants}, control the masses and the numbers of Birman--Schwinger eigenvalues at or above one in the profiles.
\begin{revision}
We first state the profile decomposition, the compactness criterion and
the optimizer system. We then give the rank-two classification and the
strict monotonicity and uniform spectral-gap consequences.
Only non-strict monotonicity from Lemma~\ref{lem:critical-constants} is
used in proving the profile decomposition; strict monotonicity is derived
afterwards in Corollary~\ref{cor:strict-mass-gap}.
\end{revision}

\begin{theorem}\label{thm:seven}
		Let \(N\in\N\), assume \eqref{eq:critical-regime}, let \(q_c\) be given by \eqref{eq:critical-exponents}, and suppose that \(V_n\geq0\) satisfies
		\begin{equation*}
			\mu_N(V_n)=1,
			\qquad
			\int_{\R^d}V_n^{q_c}\rightarrow
			\frac{N}{\ell_{0,d,s}^{(N)}}.
		\end{equation*}
		Then, after passing to a subsequence, there exist \(J\geq1\), integers \(N_1,\ldots,N_J\geq1\), critical optimizers \(V^{(j)}\neq0\), scales \(t_n^{(j)}>0\), and centers \(x_n^{(j)}\in\R^d\) such that
		\begin{equation*}
			\sum_{j=1}^JN_j=N,
			\qquad
			\mu_{N_j}(V^{(j)})=1,
			\qquad
			\frac{N}{\ell_{0,d,s}^{(N)}}
			=\sum_{j=1}^J\frac{N_j}{\ell_{0,d,s}^{(N_j)}},
		\end{equation*}
		\begin{equation*}
			\frac{t_n^{(j)}}{t_n^{(k)}}
			+\frac{t_n^{(k)}}{t_n^{(j)}}
			+t_n^{(j)}t_n^{(k)}|x_n^{(j)}-x_n^{(k)}|^2
			\rightarrow\infty
			\quad(j\neq k),
		\end{equation*}
		and
		\begin{equation*}
			\left\|V_n-
			\sum_{j=1}^J(t_n^{(j)})^{2s}
			V^{(j)}\bigl(t_n^{(j)}(\,\cdot-x_n^{(j)})\bigr)
			\right\|_{L^{q_c}}
			\rightarrow0.
		\end{equation*}
	\end{theorem}
	
	\noindent
	The proof combines the vector-valued Sobolev profile decomposition with Birman--Schwinger min--max arguments. Equality in a chain of mass inequalities identifies each profile as an optimizer, and the spectral gap in Theorem~\ref{thm:nine} gives the exact rank identity.
	
	\noindent
	For \(N=1\), the rank identity in Theorem~\ref{thm:seven} forces \(J=1\), and the strong \(L^{q_c}\) approximation in that theorem gives convergence after suitable translations and dilations.
	Hence \(\ell_{0,d,s}^{(1)}\) admits an optimizer, and every optimizing sequence \((V_n)\), normalized by \(\mu_1(V_n)=1\), is precompact in \(L^{q_c}\) modulo translations and dilations.
	\rev{For \(N\geq3\), existence and precompactness follow from
Theorem~\ref{thm:seven} under the strict subadditivity condition below.}
	
	\begin{corollary}
		\label{cor:critical-compactness}
		\rev{Assume \eqref{eq:critical-regime}, let \(N\geq3\), and suppose that}
		\begin{equation}\label{eq:critical-strict-subadditivity}
			\frac{N}{\ell_{0,d,s}^{(N)}}
			<
			\frac{K}{\ell_{0,d,s}^{(K)}}
			+
			\frac{N-K}{\ell_{0,d,s}^{(N-K)}}
			\qquad (K=1,\ldots,N-1).
		\end{equation}
		Then every optimizing sequence \((V_n)\) for \(\ell_{0,d,s}^{(N)}\), normalized by \(\mu_N(V_n)=1\), is precompact in \(L^{q_c}\) modulo translations and dilations and, \rev{after passing to a subsequence and applying suitable translations
and dilations depending on \(n\), converges to an optimizer.}
		In particular, \eqref{eq:critical-strict-subadditivity} holds if
		\begin{equation*}
			\ell_{0,d,s}^{(N)}
			>
			\max_{1\leq k<N}\ell_{0,d,s}^{(k)}.
		\end{equation*}
	\end{corollary}
	
	\begin{unnumberedremark}
		\begin{revision}
Theorem~\ref{thm:rank-two} gives \(a_2=2a_1\), so strict subadditivity
fails at rank two and the rank-two infimum is not attained.
For \(N\geq3\), Corollary~\ref{cor:critical-compactness} remains a
conditional existence result; we do not verify its hypothesis for every
such rank.
\end{revision}
	\end{unnumberedremark}
	
	\begin{theorem}\label{thm:nine}
		Assume \eqref{eq:critical-regime}, let \(N\in\N\), and suppose that \(V_*\geq0\) is an optimizer for \(\ell_{0,d,s}^{(N)}\) and is normalized by \(\mu_N(V_*)=1\).
		Then \(\mu_{N+1}(V_*)<\mu_N(V_*)=1\).
		Moreover, there exist an integer
		\begin{equation*}
			1\leq r\leq\dim\ker(\mathcal T_{V_*}-I)
		\end{equation*}
		and nonzero, pairwise \(\dot H^s\)-orthogonal functions \(f_1,\ldots,f_r\in\ker(\mathcal T_{V_*}-I)\subset\dot H^s(\R^d)\) such that
		\begin{equation}\label{eq:intro-critical-euler}
			(-\Delta)^sf_j=V_*f_j,
			\qquad
			\sum_{j=1}^r\|f_j\|_{\dot H^s}^2=1,
			\qquad
			V_*=\left(\frac{N}{\ell_{0,d,s}^{(N)}}
			\sum_{j=1}^r|f_j|^2\right)^{2s/(d-2s)}.
		\end{equation}
        The differential equation holds weakly, and the formula for the potential holds almost everywhere.
        The functions \(f_j\) and \(V_*\) have bounded continuous representatives.
        Put \(\kappa=4s/(d-2s)\).
        For every \(0<\alpha<\min\{2s,1\}\),
        \[
        f_j\in C^{0,\alpha}_{\mathrm{loc}}(\R^d),\qquad
        V_*\in C^{0,\alpha\min\{1,\kappa\}}_{\mathrm{loc}}(\R^d).
        \]
        Moreover, there are \(b_1,\ldots,b_r\in\mathbb C\) such that, uniformly in direction as \(|x|\to\infty\),
        \begin{equation}\label{eq:intro-critical-decay}
        \begin{aligned}
        f_j(x)&=|x|^{-(d-2s)}\bigl(b_j+o(1)\bigr),\\
        V_*(x)&=|x|^{-4s}\bigl(c_\infty+o(1)\bigr),\\
        c_\infty&=\left(\frac{N}{\ell_{0,d,s}^{(N)}}
        \sum_{j=1}^r|b_j|^2\right)^{2s/(d-2s)}\geq0.
        \end{aligned}
        \end{equation}

	Moreover, there exists an optimizer \(W\geq0\) for the same constant, normalized by \(\mu_N(W)=1\), such that
        \[
        \lim_{|x|\to\infty}|x|^{4s}W(x)>0.
        \]
        It can be chosen as the Kelvin transform of a translate of \(V_*\).
    \end{theorem}
	
	\noindent
	One-sided perturbations and approximate identities show that the threshold eigenvalue cluster ends at the \(N\)-th eigenvalue. Separation of the convex hull of the densities \(|g|^2\), where \(g\) is a \(\dot H^s\)-normalized eigenfunction with Birman--Schwinger eigenvalue one, then gives \eqref{eq:intro-critical-euler}.
	
	\begin{revision}
\Needspace{14\baselineskip}
\begin{theorem}\label{thm:rank-two}
Assume \eqref{eq:critical-regime}. Then
\[
 a_2=2a_1=2S_{d,s}^{q_c},\qquad
 \ell_{0,d,s}^{(2)}=\ell_{0,d,s}^{(1)},
\]
and neither the infimum defining \(a_2\) nor the supremum defining
\(\ell_{0,d,s}^{(2)}\) is attained.
If \(\mu_2(V_n)=1\) and \(\int V_n^{q_c}\to a_2\), then, after
passing to a subsequence,
\[
 \left\|V_n-\sum_{j=1}^2(t_n^{(j)})^{2s}
 Q_{d,s}\bigl(t_n^{(j)}(\,\cdot-x_n^{(j)})\bigr)\right\|_{q_c}\to0,
\]
where the two parameter families satisfy the orthogonality condition in
Theorem~\ref{thm:seven}. In particular, no normalized rank-two optimizing
sequence is precompact modulo translations and dilations.
\end{theorem}

\begin{corollary}
\label{cor:strict-mass-gap}
Assume \eqref{eq:critical-regime}. For every \(N\geq1\),
\[
 a_N<a_{N+1}.
\]
Every normalized optimizing sequence for \(\ell_{0,d,s}^{(N)}\) satisfies
\[
 \limsup_{n\to\infty}\mu_{N+1}(V_n)
 \leq\left(\frac{a_N}{a_{N+1}}\right)^{1/q_c}<1.
\]
The same bound holds for \(\mu_{N+1}(V_*)\) at every normalized optimizer
\(V_*\) of rank \(N\).
\end{corollary}

Strict increase of \(a_N=N/\ell_{0,d,s}^{(N)}\) is distinct from strict
subadditivity and does not imply monotonicity of
\(\ell_{0,d,s}^{(N)}\). In particular, the equality \(a_2=2a_1\) is
consistent with \(a_2>a_1\).
\end{revision}

The paper is organized as follows.
Section~\ref{sec:basic-results} collects the analytic preliminaries.
Section~\ref{sec:critical-euler} establishes the spectral gap and Euler--Lagrange system,
regularity and polynomial decay.
Section~\ref{sec:critical-profile} proves the critical profile decomposition\rev{,
the compactness criterion, and the rank-two and strict-monotonicity
consequences.}

\section{Preliminaries}\label{sec:basic-results}
	
	We first recall the standard analytic and perturbative tools used in the proofs.
	We then collect the spectral and compactness statements needed for the critical variational problem.
	
	For completeness, the Fourier definition of the fractional Laplacian has the singular-integral representation
\[
    (-\Delta)^sf(x)=c_{d,s}^{\mathrm{FL}}\,\mathrm{p.v.}\int_{\R^d}
    \frac{f(x)-f(y)}{|x-y|^{d+2s}}\,dy,
    \qquad f\in\mathcal S(\R^d).
\]

The fractional Sobolev and Hardy--Littlewood--Sobolev inequalities take the following forms. If \(d>2s\), then
		\begin{equation*}
			\|u\|_{L^{2_s^*}}^2
			\leq
			C_{d,s}\normHs u^2.
		\end{equation*}
		Moreover, the Riesz potential
		\begin{equation*}
			I_{2s}f(x)
			=c_{d,s}^{\mathrm{R}}\int_{\R^d}\frac{f(y)}{|x-y|^{d-2s}}\,dy
		\end{equation*}
		maps \(L^{2d/(d+2s)}\) continuously into \(L^{2d/(d-2s)}\).
Here \(c_{d,s}^{\mathrm{FL}}\) and \(c_{d,s}^{\mathrm{R}}\) denote the respective normalization constants for the Fourier multipliers \(|\xi|^{2s}\) and \(|\xi|^{-2s}\).

\begin{revision}
For real \(u,v\in\dot H^s(\R^d)\), the associated bilinear form is
\[
 \langle u,v\rangle_{\dot H^s}
 =\frac{c_{d,s}^{\mathrm{FL}}}{2}
 \iint_{\R^d\times\R^d}
 \frac{(u(x)-u(y))(v(x)-v(y))}{|x-y|^{d+2s}}\,dx\,dy.
\]
The identity extends from smooth functions by density.
Since \(0<s<1\), \(|u|,u^+,u^-\in\dot H^s\), where
\(u^+=\max\{u,0\}\) and \(u^-=\max\{-u,0\}\), and
\(\||u|\|_{\dot H^s}\leq\|u\|_{\dot H^s}\).
If \(u^+\) and \(u^-\) are both nonzero, then
\begin{equation}\label{eq:nodal-cross-term}
 \langle u^+,u^-\rangle_{\dot H^s}
 =-c_{d,s}^{\mathrm{FL}}
 \iint_{\R^d\times\R^d}
 \frac{u^+(x)u^-(y)}{|x-y|^{d+2s}}\,dx\,dy<0.
\end{equation}
The integral is finite by the bilinear-form estimate and strictly positive
before the minus sign because both factors are nontrivial and nonnegative.
\end{revision}

We also recall the first-order perturbation formula for an isolated eigenvalue of finite multiplicity; see \cite{Kato}.
For an analytic self-adjoint form family \(A_\varepsilon=A+\varepsilon B\) with a common form domain, let \(\lambda\) be an isolated eigenvalue of \(A\) of multiplicity \(r\), and let \(e_1,\ldots,e_r\) be an orthonormal basis of its eigenspace.
The first-order coefficients of the eigenvalue branches are the eigenvalues of the Hermitian matrix \((b[e_i,e_j])_{i,j=1}^r\), where \(b\) is the perturbing form.
For the critical operator \(\cT_{V+\varepsilon W}\), the perturbation matrix is
\begin{equation}\label{eq:cluster-perturbation-matrices}
\left(\int_{\R^d}Wg_i\overline{g_j}\right)_{i,j=1}^r,
\end{equation}
where \((g_i)\) is an orthonormal basis of the relevant isolated eigenspace
in \(\dot H^s\).
For an increasing ordering within a cluster, if the matrix eigenvalues are \(\beta_1\leq\cdots\leq\beta_r\), the right derivative of the \(k\)-th level is \(\beta_k\), whereas its left derivative is \(\beta_{r-k+1}\).
For a decreasing ordering, writing \(\alpha_1\geq\cdots\geq\alpha_r\), the corresponding derivatives are \(\alpha_k\) and \(\alpha_{r-k+1}\).
Thus the ordered levels need not be two-sided differentiable at a multiple eigenvalue.

Write \(N_{>1}(K)\) for the number of eigenvalues of \(K\) strictly larger than one, counted with multiplicity.
The following form of the Birman--Schwinger principle will be used throughout.

\begin{proposition}\label{prop:birman-schwinger}
\begin{enumerate}[label=\textup{(\roman*)},leftmargin=*]
\item 
If \(d>2s\) and \(V\in L^{q_c}(\R^d)\) is real, define \(\cT_V\) by the same potential form as in the nonnegative case. Then \(\cT_V\) is compact and self-adjoint on \(\dot H^s\), with
\[
    \|\cT_V\|\leq C_{d,s}\|V\|_{q_c}.
\]
It is nonnegative when \(V\geq0\).

\item 
Let \(0\leq V\in L^{q_c}(\R^d)\).
For every \(E>0\),
\[
    K_{s,V}(E):=V^{1/2}\bigl((-\Delta)^s+E\bigr)^{-1}V^{1/2}
\]
is compact, self-adjoint and nonnegative on \(L^2\), and
\[
    \#\{j:\lambda_j(H_{s,V})<-E\}=N_{>1}(K_{s,V}(E)).
\]

\item 
If \(d>2s\) and \(0\leq V\in L^{q_c}\), then \(K_{s,V}(0)\) from \eqref{eq:critical-birman-schwinger} is compact, self-adjoint and nonnegative on \(L^2(\R^d)\).
Its nonzero spectrum agrees with that of \(\cT_V\), including multiplicities, and
\[
    N_-(H_{s,V})=N_{>1}(K_{s,V}(0))=N_{>1}(\cT_V).
\]
\end{enumerate}
\end{proposition}

\begin{proof}
\emph{(i)} H\"older and Sobolev give
\[
    \left|\int V f\overline g\right|
    \leq C_{d,s}\|V\|_{q_c}\normHs f\normHs g.
\]
This defines \(\cT_V\), gives its norm bound, and proves self-adjointness and, for \(V\geq0\), positivity.
For bounded compactly supported \(V\), local Rellich compactness gives compactness of \(\cT_V\).
Approximation in \(L^{q_c}\), together with the same norm estimate for \(V-V_k\), proves the general case.

\emph{(ii)} Put \(A=(-\Delta)^s\) and \(B_E=V^{1/2}(A+E)^{-1/2}\).
H\"older and the critical Sobolev inequality imply
\(\|B_Ef\|_2^2\leq C_E\|V\|_{q_c}\|f\|_2^2\).
Local Rellich compactness proves compactness for bounded compactly supported potentials; approximation extends it to all admissible \(V\), since
\(|\sqrt V-\sqrt W|^2\leq|V-W|\).
Thus \(K_{s,V}(E)=B_EB_E^*\) is a well-defined compact nonnegative operator.
The bijection \(u\mapsto f=(A+E)^{1/2}u\) from \(H^s\) onto \(L^2\) gives
\[
    \normHs u^2+E\|u\|_2^2-\int V|u|^2
    =\langle f,(I-B_E^*B_E)f\rangle_{L^2}.
\]
The variational principle identifies the maximal negative dimensions of these forms.
Since \(B_E^*B_E\) and \(B_EB_E^*\) have the same nonzero eigenvalues with multiplicities, the counting identity follows.

\emph{(iii)} The same Sobolev--Rellich approximation argument makes
\(J:\dot H^s\to L^2\), \(Ju=V^{1/2}u\), compact.
With the adjoint taken between these Hilbert spaces,
\[
    \cT_V=J^*J,\qquad K_{s,V}(0)=JJ^*.
\]
These factorizations define the zero-energy product and identify the nonzero spectra.
Let \(U:\dot H^s\to L^2\), \(Uu=A^{1/2}u\), be unitary, and put \(C=JU^{-1}\).
Then
\[
    K_{s,V}(E)=C\,A(A+E)^{-1}C^*,\qquad K_{s,V}(0)=CC^*.
\]
As \(E\searrow0\), the contractions \(A(A+E)^{-1}\) converge strongly to \(I\), since \(\ker_{L^2}A=\{0\}\); compactness of \(C\) therefore gives norm convergence \(K_{s,V}(E)\to K_{s,V}(0)\).
The operator inequality \(0\leq K_{s,V}(E)\leq K_{s,V}(0)\) gives \(N_{>1}(K_{s,V}(E))\leq N_{>1}(K_{s,V}(0))\). Compactness makes the latter count finite, and norm convergence ensures that all these eigenvalues remain strictly above one for sufficiently small \(E>0\). Hence the two counts agree for all sufficiently small \(E\).
Applying \textup{(ii)} to the critical potential and using
\(N_-(H_{s,V}+E)\nearrow N_-(H_{s,V})\) completes the proof.
\end{proof}

To verify \eqref{eq:critical-truncated-relation}, apply the preceding counting identity.
	Indeed, if \(r=\min\{N,N_-(H_{s,V})\}\geq1\), then \(\mu_r(V)>1\), which gives the upper bound in \eqref{eq:critical-truncated-relation}, while the lower bound follows by scaling \(V\) by \(c>\mu_k(V)^{-1}\) and letting \(c\searrow\mu_k(V)^{-1}\).

The preceding spectral tools yield the basic properties of the critical constants.

\begin{lemma}\label{lem:critical-constants}
Assume \eqref{eq:critical-regime}, and let \(a_N\) be defined by \eqref{eq:critical-minimization}.
Then \(0<a_N<\infty\) for every \(N\geq1\). Moreover,
\begin{equation}\label{eq:critical-rank-monotonicity}
    a_N\leq a_{N+1}\qquad(N\geq1),
\end{equation}
and
\begin{equation}\label{eq:critical-subadditivity}
    a_N\leq a_K+a_{N-K}\qquad(1\leq K\leq N-1).
\end{equation}
In particular, \(\ell_{0,d,s}^{(N)}\geq\ell_{0,d,s}^{(1)}\), and
\begin{equation}\label{eq:critical-limit}
    \lim_{N\to\infty}\ell_{0,d,s}^{(N)}
    =\sup_{N\geq1}\ell_{0,d,s}^{(N)}
    =L_{0,d,s}.
\end{equation}
\end{lemma}

\begin{proof}
		We first establish positivity and finiteness. For every admissible \(V\),
the Sobolev estimate in Proposition~\ref{prop:birman-schwinger} gives
\[
 1\leq\mu_N(V)\leq\|\cT_V\|\leq C_{d,s}\|V\|_{q_c},
\]
and hence \(a_N\geq C_{d,s}^{-q_c}>0\).
To show that the admissible set is nonempty, choose a ball \(B\),
a nonnegative \(V_0\in C_c^\infty(\R^d)\) with \(V_0\geq1\) on \(B\),
and an \(N\)-dimensional subspace \(F\subset C_c^\infty(B)\).
The continuous quadratic form \(f\mapsto\int V_0|f|^2\) is strictly
positive on the \(\dot H^s\)-unit sphere of \(F\). By finite-dimensional
compactness,
\[
 c_F:=\min_{\substack{f\in F\\\|f\|_{\dot H^s}=1}}
       \int V_0|f|^2>0.
\]
The min--max principle gives \(\mu_N(V_0)\geq c_F\), so
\(c_F^{-1}V_0\) is admissible and
\(a_N\leq c_F^{-q_c}\int V_0^{q_c}<\infty\).

By \eqref{eq:critical-minimization} and \(\mu_{N+1}(V)\leq\mu_N(V)\), the admissible sets decrease with \(N\), which proves \eqref{eq:critical-rank-monotonicity}.
		
		For \eqref{eq:critical-subadditivity}, choose near-minimizers \(V\) and \(W\) for \(a_K\) and \(a_{N-K}\), respectively, and let \(W_R(x):=W(x-Re_1)\).
		On \(\dot H^s\), \(\mathcal T_{V+W_R}=\mathcal T_V+U_R\mathcal T_WU_R^*\), where \(U_R\) is the unitary translation.
		Translations converge weakly to zero, while \(\mathcal T_V\) and \(\mathcal T_W\) are compact.
		Define
		\begin{equation*}
			C_R(f,g):=\mathcal T_V^{1/2}f
			+U_R\mathcal T_W^{1/2}g.
		\end{equation*}
		Then \(C_RC_R^*=\mathcal T_{V+W_R}\), whereas the diagonal blocks of \(C_R^*C_R\) are \(\mathcal T_V\) and \(\mathcal T_W\), and its off-diagonal blocks tend to zero in operator norm by compactness and weak convergence of translations.
		Thus \(C_R^*C_R\to\mathcal T_V\oplus\mathcal T_W\) in norm.
		The nonzero spectra of \(C_RC_R^*\) and \(C_R^*C_R\) agree with multiplicity, which proves convergence to the ordered spectral union.
		It follows that \(\mu_N(V+W_R)\geq1-o(1)\).
		Multiplying \(V+W_R\) by a scalar tending to one makes it admissible for \(a_N\).
		Approximating \(V\) and \(W\) in \(L^{q_c}\) by bounded compactly supported functions, whose translated supports are eventually disjoint, gives
		\begin{equation*}
			\int_{\R^d}(V+W_R)^{q_c}
			=
			\int_{\R^d}V^{q_c}+\int_{\R^d}W^{q_c}+o(1).
		\end{equation*}
		Letting first \(R\to\infty\) and then the near-minimizer errors tend to zero proves \eqref{eq:critical-subadditivity}.
		Iteration yields \(a_N\leq Na_1\), hence \(\ell_{0,d,s}^{(N)}\geq \ell_{0,d,s}^{(1)}\).
		
		By Fekete's lemma and \eqref{eq:critical-subadditivity},
		\begin{equation*}
			\lim_{N\to\infty}\frac{a_N}{N}
			=
			\inf_{N\geq1}\frac{a_N}{N},
		\end{equation*}
		so the first equality in \eqref{eq:critical-limit} follows after taking reciprocals.
		It remains to identify the supremum.
		Given \(V\) with \(\mu_N(V)>0\), apply the Birman--Schwinger principle to \(cV\), where \(c>\mu_N(V)^{-1}\).
		Then \(N_-(H_{s,cV})\geq N\), and therefore
		\begin{equation*}
			L_{0,d,s}
			\geq
			\frac{N}{c^{q_c}\int V^{q_c}}.
		\end{equation*}
		Letting \(c\searrow\mu_N(V)^{-1}\) and optimizing in \(V\) gives \(L_{0,d,s}\geq\ell_{0,d,s}^{(N)}\).
		Conversely, if \(M=N_-(H_{s,V})>0\), then \(\mu_M(V)>1\), whence
		\begin{equation*}
			\ell_{0,d,s}^{(M)}
			\geq
			\frac{M\mu_M(V)^{q_c}}{\int V^{q_c}}
			\geq
			\frac{M}{\int V^{q_c}}.
		\end{equation*}
		Taking the supremum over \(V\) proves the reverse inequality and completes \eqref{eq:critical-limit}.
		This follows the argument of \cite[Lemma~6]{FGL}, with separated-well spectral convergence supplied by the fractional compact-operator decoupling above.
	\end{proof}

Concentration may occur through changes of position and scale.
The following vector profile decomposition tracks translations and dilations simultaneously for a finite family of functions; it will be applied to Birman--Schwinger eigenfunctions in Section~\ref{sec:critical-profile}.

\begin{lemma}\label{lem:vector-profile}
Let \(M\geq1\) be fixed, \(0<s<\min\{1,d/2\}\), and let
\(F_n=(f_{1,n},\ldots,f_{M,n})\) be bounded in
\(\mathcal H:=\dot H^s(\R^d;\C^M)\), with
\(\|F\|_{\mathcal H}^2=\sum_{j=1}^M\normHs{f_j}^2\).
After passing to a subsequence, there are \(J\in\{0,1,\ldots\}\cup\{\infty\}\) nonzero profiles
\(F^{(k)}\in\mathcal H\) and parameters \(t_n^{(k)}>0\), \(x_n^{(k)}\in\R^d\), \(1\leq k\leq J\), such that
\[
    \frac{t_n^{(k)}}{t_n^{(\ell)}}+
    \frac{t_n^{(\ell)}}{t_n^{(k)}}+
    t_n^{(k)}t_n^{(\ell)}|x_n^{(k)}-x_n^{(\ell)}|^2
    \rightarrow\infty\qquad(k\neq\ell).
\]
Write \(T_n^{(k)}f=(t_n^{(k)})^{(d-2s)/2}
f(t_n^{(k)}(\,\cdot-x_n^{(k)}))\), acting on each component.
For each fixed \(L\geq1\), put \(L_J=\min\{L,J\}\). Then
\[
    f_{j,n}=\sum_{k=1}^{L_J}T_n^{(k)}f_j^{(k)}+r_{j,n}^{(L)},
    \qquad 1\leq j\leq M,
\]
and, with \(R_n^{(L)}=(r_{1,n}^{(L)},\ldots,r_{M,n}^{(L)})\),
\begin{equation}\label{eq:vector-remainder-orthogonality}
    (T_n^{(k)})^{-1}R_n^{(L)}\rightharpoonup0
    \quad\text{in }\mathcal H\qquad(1\leq k\leq L_J).
\end{equation}
For every fixed \(c=(c_1,\ldots,c_M)\in\C^M\) and \(L\),
\begin{align}\label{eq:vector-linear-combination-decoupling}
    \left\|\sum_{j=1}^M c_jf_{j,n}\right\|_{\dot H^s}^2
    &=\sum_{k=1}^{L_J}
      \left\|\sum_{j=1}^M c_jf_j^{(k)}\right\|_{\dot H^s}^2
      +\left\|\sum_{j=1}^M c_jr_{j,n}^{(L)}\right\|_{\dot H^s}^2
      +o_n(1).
\end{align}
Moreover,
\[
    \lim_{L\to\infty}\limsup_{n\to\infty}
    \|r_{j,n}^{(L)}\|_{L^{2_s^*}}=0\qquad(1\leq j\leq M).
\]
If \(J<\infty\), the remainder is unchanged for \(L\geq J\); if \(J=0\), the sums are empty.
\end{lemma}

\begin{proof}
We use the extraction in \cite[Theorem~1.4 and Section~5]{PalatucciPisante}.
The maps \(T_{t,x}f=t^{(d-2s)/2}f(t(\,\cdot-x))\) act unitarily on
\(\dot H^s\) and diagonally on \(\mathcal H\).
We make quantitative the concentration criterion of
\cite[Lemma~4.1]{PalatucciPisante}. Put \(P=2_s^*\), \(\theta=2/P\), and
\[
 \mathfrak m(u)=\sup_{y\in\R^d,\,R>0}
 \left(R^{-2s}\int_{B_R(y)}|u|^2\right)^{1/2}.
\]
The refined Sobolev inequality \cite[Theorem~1.1, with \(r=2\)]{PalatucciPisante}
gives \(\|u\|_P\leq C\|u\|_{\dot H^s}^{\theta}
\mathfrak m(u)^{1-\theta}\) for real \(u\).
If \(\|u_n\|_{\dot H^s}\leq B\) and
\(\limsup_n\|u_n\|_P\geq\varepsilon>0\), pass to a subsequence with
\(\|u_n\|_P\geq\varepsilon/2\). Then
\[
 \mathfrak m(u_n)\geq
 \delta:=\left(\frac{\varepsilon}{2CB^\theta}\right)^{1/(1-\theta)}>0.
\]
Choose \(y_n,R_n\) so that the defining squared expression is at least
\(\delta^2/2\), and set
\[w_n(x)=R_n^{(d-2s)/2}u_n(y_n+R_nx).\]
Weak compactness and local Rellich compactness yield
\(w_n\rightharpoonup w\) in \(\dot H^s\) and \(w_n\to w\) in \(L^2(B_1)\).
Thus \(\int_{B_1}|w|^2\geq\delta^2/2\). By H\"older and Sobolev,
\(\|w\|_{\dot H^s}\geq c_{d,s}\delta\).
This is the required bound \(c(B,\varepsilon)>0\), uniform across extraction steps.
For complex sequences, select a real or imaginary part with \(L^P\)-norm
at least \(\varepsilon/4\), then extract the full complex weak limit.

For a bounded vector sequence \(R_n\), let \(\eta(R)\) be the supremum of the \(\mathcal H\)-norms of all weak limits of \(T_n^{-1}R_n\), allowing subsequences and arbitrary translation--dilation parameters.
If \(\eta(R)>0\), choose a weak limit \(F\) with \(\|F\|_{\mathcal H}\geq\eta(R)/2\), and subtract \(T_nF\) from every component.
Weak convergence and unitarity give
\[
    \|R_n\|_{\mathcal H}^2
    =\|F\|_{\mathcal H}^2+\|R_n-T_nF\|_{\mathcal H}^2+o_n(1),
    \qquad T_n^{-1}(R_n-T_nF)\rightharpoonup0.
\]
Iterate this selection, stopping when \(\eta=0\).
A new parameter sequence must diverge from every previous one in the sense stated in the lemma: otherwise the relative transformations have a strongly convergent subsequence, and weak vanishing in that earlier frame would force the proposed new profile to be zero.
Conversely, divergent parameters give weakly vanishing relative transformations, so subtracting a new profile preserves weak vanishing in all earlier frames.
The usual diagonal subsequence construction in \cite[Section~5]{PalatucciPisante} gives a single subsequence with these properties for every fixed extraction index.
In particular, it yields \eqref{eq:vector-remainder-orthogonality} and
\[
    \sum_{k=1}^J\|F^{(k)}\|_{\mathcal H}^2
    \leq\limsup_n\|F_n\|_{\mathcal H}^2.
\]

To obtain the stronger componentwise identity, parameter separation gives
\[
    \langle T_n^{(k)}f_i^{(k)},T_n^{(\ell)}f_j^{(\ell)}\rangle_{\dot H^s}
    \rightarrow0\qquad(k\neq\ell)
\]
for every \(i,j\).
Equation~\eqref{eq:vector-remainder-orthogonality} also gives
\[
    \langle T_n^{(k)}f_i^{(k)},r_{j,n}^{(L)}\rangle_{\dot H^s}
    \rightarrow0\qquad(k\leq L_J).
\]
Expanding the squared norm of any fixed linear combination and using unitarity proves \eqref{eq:vector-linear-combination-decoupling} on the same subsequence for all \(c\).

Finally, the energy bound and the choice of profiles within a factor two of the maximal defect imply that the defects of the successive remainders tend to zero; if extraction stops, the final defect is zero.
The remainders have a uniform limiting upper \(\mathcal H\)-bound by the energy identity.
If the asserted \(L^{2_s^*}\)-remainder limit failed, then, since \(M\) is fixed, some component along a sequence of extraction indices would have limiting upper norm at least \(\varepsilon>0\).
The scalar concentration criterion would give parameters and a weak limit of that component with norm at least \(c(B,\varepsilon)>0\), uniformly in the extraction index.
Weak compactness of the other finitely many components then gives a vector weak limit with at least the same norm, contradicting the vanishing defect.
This proves the remainder assertion without applying a scalar Sobolev estimate to the vector modulus.
\end{proof}

\section{Euler--Lagrange equations and decay of critical optimizers}\label{sec:critical-euler}
\begin{revision}
We adapt the critical perturbation argument of \cite[Section~8]{FGL}
to \(\mathcal T_V\) on \(\dot H^s\). We first establish the threshold
spectral gap and the Euler--Lagrange system. The subsequent regularity
and decay argument uses the fractional Riesz kernel and its Kelvin
covariance.
\end{revision}

At the critical threshold, \(\cT_Vu=u\) is equivalent to the weak equation \((-\Delta)^su=Vu\) in \(\dot H^s\).
The factorization \(\cT_V=J^*J\), \(K_{s,V}(0)=JJ^*\), with \(Ju=V^{1/2}u\), identifies the eigenspaces at one: \(u\mapsto Ju\) and \(\phi\mapsto J^*\phi\) are inverse on these eigenspaces.
A nonzero such solution in \(L^2\) is a zero-energy eigenfunction; one in \(\dot H^s\setminus L^2\) is called a zero-energy resonance here.
Neither contributes to the strict negative-spectrum count in Proposition~\ref{prop:birman-schwinger}.
	
	The spectral gap will later be used to obtain the exact rank identity in the profile decomposition. The argument below treats a multiple threshold eigenvalue using one-sided perturbations.
	
	For later perturbations by a real \(U\in L^{q_c}(\R^d)\), not necessarily nonnegative, we extend the same notation to the \(N\)-th variational level of the compact self-adjoint operator \(\cT_U\), namely
	\begin{equation*}
		\mu_N(U):=
		\sup_{\substack{F\subset\dot H^s(\R^d)\\ \dim F=N}}
		\ \inf_{\substack{f\in F\\ \normHs f=1}}
		\ipHs f{\cT_Uf}.
	\end{equation*}
	This convention agrees with the preceding eigenvalue notation when \(U\geq0\); below it is used only where the \(N\)-th variational level is positive, as holds in a sufficiently small neighborhood of the reference potential.

	\begin{proof}[Proof of the spectral and Euler assertions in Theorem~\ref{thm:nine}]
		Let
		\begin{equation*}
			E:=\ker(\cT_{V_*}-I),
			\qquad
			m_*:=\dim E,
		\end{equation*}
		and let \(g_1,\ldots,g_{m_*}\) be a \(\dot H^s\)-orthonormal basis of \(E\).
		Let \(r_>\) be the number of eigenvalues strictly larger than one, and put \(K:=N-r_>\).
		Thus the \(N\)-th eigenvalue is the \(K\)-th selected eigenvalue in the threshold block.
		
		\medskip
		\noindent\emph{Step 1: one-sided perturbation inequalities.}
		For real \(W\in L^{q_c}\), define
		\begin{equation*}
			A_W
			:=
			\left(\int Wg_i\overline{g_j}\right)_{i,j=1}^{m_*},
			\qquad
			c(W)
			:=
			\frac{\displaystyle\int V_*^{q_c-1}W}
			{\displaystyle\int V_*^{q_c}}.
		\end{equation*}
		Although \(V_*+\varepsilon W\) need not be nonnegative, the form operator \(\cT_{V_*+\varepsilon W}\) is well defined and compact.
		For sufficiently small \(|\varepsilon|\), its \(N\)-th eigenvalue remains positive.
		Consider
		\begin{equation*}
			\Phi_W(\varepsilon):=
			\frac{\mu_N(V_*+\varepsilon W)^{q_c}}
			{\displaystyle\int_{\R^d}(V_*+\varepsilon W)_+^{q_c}}.
		\end{equation*}
		The form inequality \(\cT_{V_*+\varepsilon W}\leq \cT_{(V_*+\varepsilon W)_+}\), followed by the variational definition of \(\ell_{0,d,s}^{(N)}\), shows that
		\begin{equation*}
			\Phi_W(\varepsilon)\leq\frac{\ell_{0,d,s}^{(N)}}N
			=\Phi_W(0).
		\end{equation*}
		Thus zero is a two-sided maximum for every real \(W\in L^{q_c}\), without any support restriction.
		Since \(q_c>1\), the positive-part map satisfies
		\begin{equation*}
			\int(V_*+\varepsilon W)_+^{q_c}
			=\int V_*^{q_c}
			+q_c\varepsilon\int V_*^{q_c-1}W+o(|\varepsilon|).
		\end{equation*}
		Let \(\alpha_1(W)\geq\cdots\geq\alpha_{m_*}(W)\) denote the eigenvalues of \(A_W\).
		The one-sided perturbation formulas following \eqref{eq:cluster-perturbation-matrices} give
		\begin{equation*}
			\mu_N(V_*+\varepsilon W)=
			\begin{cases}
				1+\varepsilon\alpha_K(W)+o(\varepsilon),&\varepsilon\searrow0,\\
				1+\varepsilon\alpha_{m_*-K+1}(W)+o(|\varepsilon|),&\varepsilon\nearrow0.
			\end{cases}
		\end{equation*}
		\rev{Expanding \(\Phi_W\), and using its nonpositive right derivative and
nonnegative left derivative at zero, gives}
		\begin{equation*}
			\alpha_K(W)
			\leq c(W)
			\leq
			\alpha_{m_*-K+1}(W),
		\end{equation*}
		for every real \(W\in L^{q_c}\).
		The use of the positive part in \(\Phi_W\) is what handles the obstacle \(V\geq0\), including perturbations on \(\{V_*=0\}\).
		
		\medskip
		\noindent\emph{Step 2: the gap below the threshold eigenvalue cluster.}
		Assume \(K<m_*\).
		Choose a common Lebesgue point \(x\) of all products \(g_i\overline{g_j}\) and of \(V_*^{q_c-1}\) such that \(V_*(x)>0\), \rev{and set \(W_{\delta,x}(y)=\delta^{-d}\rho((y-x)/\delta)\),
where \(0\leq\rho\in C_c^\infty\) and \(\int\rho=1\).
For each fixed \(\delta>0\), Step~1 applies to the direction
\(W_{\delta,x}\in L^{q_c}\); we then let \(\delta\searrow0\).}
		Then
		\begin{equation*}
			A_{W_{\delta,x}}
			\rightarrow
			G(x)G(x)^*,
			\qquad
			G(x):=(g_1(x),\ldots,g_{m_*}(x))^T.
		\end{equation*}
		This matrix has eigenvalues \(|G(x)|^2,0,\ldots,0\).
		If \(K=1<m_*\), the limiting inequality becomes \(|G(x)|^2\leq c_x\leq0\), hence \(G(x)=0\) and \(c_x=0\).
		If \(2\leq K<m_*\), both selected eigenvalues vanish and again \(c_x=0\).
		Since \(c_x=V_*^{q_c-1}(x)/\int V_*^{q_c}>0\), either conclusion is a contradiction.
		Hence \(K=m_*\), so \(N=r_>+m_*\): the threshold eigenvalue cluster ends at the \(N\)-th eigenvalue, and \(\mu_{N+1}(V_*)<1=\mu_N(V_*)\).
		
		\medskip
		\noindent\emph{Step 3: convex Euler multiplier.}
		After Step~2, \(\mu_N\) is the smallest eigenvalue in the threshold block.
We use the preceding directional inequalities directly.
Work in the real Banach space \(X=L^{q_c'}(\R^d;\R)\), and put
\[
 h_*:=\frac{V_*^{q_c-1}}{\int V_*^{q_c}},\qquad
 \mathcal D:=\operatorname{conv}
 \{\,|g|^2:g\in E,\ \|g\|_{\dot H^s}=1\,\}\subset X.
\]
Sobolev embedding places every density in \(X\). These densities lie in
the finite-dimensional real span of the real and imaginary parts of
\(g_i\overline{g_j}\). The unit sphere of \(E\) is compact, and the map
\(g\mapsto|g|^2\) is continuous into \(X\); its convex hull \(\mathcal D\)
is therefore compact as well.
For every real \(W\in L^{q_c}\), the Rayleigh principle and Step~1 give
\[
 \min_{\substack{g\in E\\\|g\|_{\dot H^s}=1}}\int W|g|^2
 =\alpha_{m_*}(W)\leq\int Wh_*\leq\alpha_1(W)
 =\max_{\substack{g\in E\\\|g\|_{\dot H^s}=1}}\int W|g|^2.
\]
If \(h_*\notin\mathcal D\), strict separation of a point from a closed
convex set in \(X\) would supply a continuous real linear functional
violating one of these inequalities. Since \(1<q_c'<\infty\), that
functional is integration against some real \(W\in L^{q_c}\).
Hence \(h_*\in\mathcal D\). Expressing this convex combination in the
basis \((g_i)\) yields a positive semidefinite matrix of trace one.
The spectral gap proved in Step~2 does not assert that this multiplier has full rank.
Consequently there is a positive semidefinite matrix \(\Gamma=(\Gamma_{ij})\) with
		\begin{equation*}
			\Tr\Gamma=1,
			\qquad
			\frac{V_*^{q_c-1}}{\int V_*^{q_c}}
			=
			\sum_{i,j=1}^{m_*}\Gamma_{ij}g_i\overline{g_j}.
		\end{equation*}
		Because \(\int V_*^{q_c}=N/\ell_{0,d,s}^{(N)}\), we obtain
		\begin{equation*}
			V_*^{q_c-1}
			=
			\frac{N}{\ell_{0,d,s}^{(N)}}
			\sum_{i,j}\Gamma_{ij}g_i\overline{g_j}.
		\end{equation*}
		After a unitary change of the \(\dot H^s\)-orthonormal basis of \(E\), diagonalize \(\Gamma\).
		Set \(r:=\operatorname{rank}\Gamma\), retain the positive eigenvalues \(\alpha_1,\ldots,\alpha_r>0\), and discard the zero eigenspaces.
		Then \(1\leq r\leq\dim E\), \(\sum_{j=1}^r\alpha_j=1\), and, with
		\begin{equation*}
			f_j:=\sqrt{\alpha_j}\,g_j\quad(j=1,\ldots,r),
			\qquad
			\sum_{j=1}^r\normHs{f_j}^2=1,
			\qquad
			\frac1{q_c-1}=\frac{2s}{d-2s},
		\end{equation*}
		the functions \(f_1,\ldots,f_r\) are nonzero and pairwise orthogonal.
		The asserted nonlinear formula follows.
		
		\medskip
		\noindent\emph{Step 4: zero-energy equation.}
		For \(g\in E\),
		\begin{equation*}
			\ipHs\varphi g
			=
			\int V_*\varphi\overline g
			\quad(\varphi\in\dot H^s),
			\qquad
			(-\Delta)^sg=V_*g
		\end{equation*}
		in the weak sense, and hence the same holds for the \(f_j\)'s.
	\end{proof}
	
	We now prove the regularity and asymptotic assertions of Theorem~\ref{thm:nine} using the Euler system.

\begin{proposition}\label{prop:critical-decay}
Let \(V_*\) and \(f_1,\ldots,f_r\) satisfy the spectral and Euler assertions of Theorem~\ref{thm:nine}.
Then its regularity and decay assertions hold.
In addition, \(f_j\in L^p(\R^d)\) for every finite \(p\geq2d/(d-2s)\), and
\[
|f_j(x)|\leq C(1+|x|)^{-(d-2s)},\qquad
V_*(x)\leq C(1+|x|)^{-4s}.
\]
The coefficient \(c_\infty\) is positive exactly when at least one \(b_j\) is nonzero.
Moreover, the conformal orbit of \(V_*\) contains an optimizer \(W\) with \(\mu_N(W)=1\) and a strictly positive leading decay coefficient.
\end{proposition}

\begin{proof}
Put
\[
q=\frac d{2s},\quad p_*=\frac{2d}{d-2s},\quad
t_0=\frac{2d}{d+2s},\quad \kappa=\frac{4s}{d-2s},\quad
a=\left(\frac{N}{\ell_{0,d,s}^{(N)}}\right)^{2s/(d-2s)}.
\]
Writing \(F=(f_1,\ldots,f_r)\), the Euler system reads
\(( -\Delta)^sf_j=V_*f_j\), with \(V_*=a|F|^\kappa\).

\medskip\noindent\emph{Step 1: the integral equation.}
Sobolev embedding gives \(F\in L^{p_*}\), and H\"older's inequality gives
\(H_j:=V_*f_j\in L^{t_0}\).
The Riesz potential \(I_{2s}=(-\Delta)^{-s}\) maps \(L^{t_0}\) into \(\dot H^s\) and into \(L^{p_*}\).
Indeed, the first mapping follows by duality from the Sobolev inequality, and the second is HLS.
Uniqueness of the weak solution in \(\dot H^s\), obtained by testing the homogeneous equation with the solution itself, therefore yields
\begin{equation}\label{eq:critical-riesz-system}
f_j(x)=I_{2s}(V_*f_j)(x)
=c_{d,s}^{\mathrm{R}}\int_{\R^d}|x-y|^{-(d-2s)}V_*(y)f_j(y)\,dy
\end{equation}
almost everywhere.

\medskip\noindent\emph{Step 2: improvement of integrability at the critical exponent.}
Suppose \(f_j\in L^{p_0}\), where \(p_0\geq p_*\), and choose a finite \(p_1>p_0\) such that
\[
\frac1{p_1}+\frac{2s}{d}\geq\frac1{p_0}.
\]
Since \(V_*\in L^q\), write \(V_*=V_b+V_\varepsilon\), where \(V_b\) is bounded and compactly supported and \(\|V_\varepsilon\|_q<\varepsilon\).
For \(Th=I_{2s}(V_\varepsilon h)\), H\"older and HLS imply
\[
\|Th\|_p\leq C_p\varepsilon\|h\|_p,\qquad p=p_0,p_1.
\]
The HLS input exponent is \(t(p)=(1/p+2s/d)^{-1}\), which lies in \((1,q)\), because \(p\geq p_*>d/(d-2s)\).
Choose \(\varepsilon\) so that both operator norms are less than \(1/2\).
As \(t(p_0),t(p_1)\leq p_0\), boundedness and compact support of \(V_b\) give
\[
g_j:=I_{2s}(V_bf_j)\in L^{p_0}\cap L^{p_1}.
\]
Equation~\eqref{eq:critical-riesz-system} is \(f_j=Tf_j+g_j\).
The Neumann series \(\sum_{n\geq0}T^ng_j\) converges in both spaces; its limits agree as distributions, and uniqueness in \(L^{p_0}\) identifies the limit with \(f_j\).
Thus \(f_j\in L^{p_1}\).
Starting at \(p_*\), a finite chain of such exponents reaches any prescribed finite \(p\geq p_*\).
This argument uses only the integral equation and the initial \(L^{p_*}\) and \(L^q\) bounds, a fact used in Step 4.

\medskip\noindent\emph{Step 3: boundedness and local continuity.}
Since \(|H_j|\leq a|F|^{\kappa+1}\), Step 2 gives \(H_j\in L^t\) for arbitrarily large finite \(t\), as well as \(H_j\in L^{t_0}\).
Choose \(t>q\) with \((\kappa+1)t\geq p_*\).
Splitting the Riesz kernel at radius one gives
\[
\|f_j\|_\infty\leq C\left(
\|H_j\|_t\big\||z|^{-(d-2s)}\big\|_{L^{t'}(B_1)}
+\|H_j\|_{t_0}\big\||z|^{-(d-2s)}\big\|_{L^{t_0'}(B_1^c)}\right)<\infty.
\]
Here \(t'\) and \(t_0'\) denote conjugate exponents; integrability at zero follows from \(t>q\), and at infinity from \((d-2s)t_0'=2d>d\).
Consequently \(V_*\) and \(H_j\) are bounded as well.

For completeness, fix a ball and choose a smooth compactly supported cutoff \(\chi\) equal to one on a larger concentric ball.
Decompose \(f_j=I_{2s}(\chi H_j)+I_{2s}((1-\chi)H_j)\).
For the first term, split the difference of the kernels into distances less than \(2|h|\) and greater than \(2|h|\).
Boundedness of \(H_j\), followed by radial integration and the mean value theorem, bounds its increment by
\[
C\begin{cases}
|h|^{2s},&2s<1,\\
|h|(1+|\log|h||),&2s=1,\\
|h|,&2s>1,
\end{cases}\qquad 0<|h|<\tfrac12.
\]
For the second term the kernel is separated from its singularity on the smaller ball.
Its gradient is integrable against \(H_j\) by H\"older with exponent \(t_0\); indeed \((d-2s+1)t_0'>d\).
This term is locally Lipschitz.
It follows that \(f_j\in C^{0,\alpha}_{\mathrm{loc}}\) for every \(0<\alpha<\min\{2s,1\}\).
The map \(z\mapsto |z|^\kappa\) is \(\min\{1,\kappa\}\)-H\"older on bounded sets, giving the asserted regularity of \(V_*\).

\medskip\noindent\emph{Step 4: regularity of the Kelvin transforms.}
For \(x\neq0\), define
\[
\widetilde f_j(x)=|x|^{-(d-2s)}f_j\left(\frac{x}{|x|^2}\right),\qquad
\widetilde V(x)=|x|^{-4s}V_*\left(\frac{x}{|x|^2}\right).
\]
Change of variables gives \(\|\widetilde f_j\|_{p_*}=\|f_j\|_{p_*}\),
\(\|\widetilde V\|_q=\|V_*\|_q\), and \(\widetilde V=a|\widetilde F|^\kappa\).
Using
\[
\left|\frac{x}{|x|^2}-\frac{z}{|z|^2}\right|
=\frac{|x-z|}{|x||z|}
\]
and the inversion Jacobian \(|z|^{-2d}\) in \eqref{eq:critical-riesz-system} yields
\begin{equation}\label{eq:critical-kelvin-integral}
\widetilde f_j=I_{2s}(\widetilde V\widetilde f_j)
\quad\hbox{almost everywhere on }\R^d.
\end{equation}
The substitutions are justified by absolute convergence almost everywhere, which follows by HLS applied to \(|V_*f_j|\) and \(|\widetilde V\widetilde f_j|\) in \(L^{t_0}\).
The single point zero does not affect this global integral identity.
Steps 2 and 3 therefore apply directly to the transformed system and give bounded continuous representatives, including at the origin.
No removable-singularity or continuity assumption is needed.

\medskip\noindent\emph{Step 5: asymptotics.}
Set \(b_j=\widetilde f_j(0)\).
Inverting the transform and using continuity at zero gives
\[
f_j(x)=|x|^{-(d-2s)}\bigl(b_j+o(1)\bigr)
\]
uniformly in direction.
Since \((d-2s)\kappa=4s\), the identity \(V_*=a|F|^\kappa\) gives \eqref{eq:intro-critical-decay}.
Boundedness on compact sets and these asymptotics imply the global upper bounds.
The formula for \(c_\infty\) proves the stated criterion for its positivity.

\medskip\noindent\emph{Step 6: an optimizer with a positive leading coefficient.}
The continuous vector \(F\) is not identically zero. Choose \(x_0\in\R^d\) with \(|F(x_0)|>0\), so that \(V_*(x_0)=a|F(x_0)|^\kappa>0\).
For \(x\neq0\), put
\begin{equation}\label{eq:critical-positive-representative}
\Phi(x)=x_0+\frac{x}{|x|^2},\qquad
W(x)=|x|^{-4s}V_*(\Phi(x)).
\end{equation}
The value at zero is immaterial for the variational problem.
Since \(|\det D\Phi(x)|=|x|^{-2d}\) and \(4sq=2d\), change of variables yields
\[
\int_{\R^d}W^q=\int_{\R^d}V_*^q.
\]
To verify preservation of the spectral normalization, define the unitary operator on \(L^2(\R^d)\)
\[
(Uh)(x)=|x|^{-d}h(\Phi(x)).
\]
The integral kernel of \(K_{s,V_*}(0)\) is
\(c_{d,s}^{\mathrm{R}}V_*(z)^{1/2}|z-w|^{-(d-2s)}V_*(w)^{1/2}\).
The identity
\[
|\Phi(x)-\Phi(y)|=\frac{|x-y|}{|x||y|}
\]
therefore gives
\begin{equation}\label{eq:critical-kelvin-conjugacy}
U K_{s,V_*}(0)U^*=K_{s,W}(0).
\end{equation}
Indeed, the transformed kernel has factors
\(|x|^{-d}|x|^{d-2s}V_*(\Phi(x))^{1/2}=W(x)^{1/2}\)
and the analogous factor in \(y\).
The kernel calculation extends to the bounded operators by HLS and density.
Thus \(\mu_j(W)=\mu_j(V_*)\) for every \(j\); in particular, \(\mu_N(W)=1\).
Together with preservation of the \(L^q\) integral, this proves that \(W\) is an optimizer for \(\ell_{0,d,s}^{(N)}\).
Finally, continuity of \(V_*\) at \(x_0\) implies, uniformly in direction,
\[
\lim_{|x|\to\infty}|x|^{4s}W(x)
=\lim_{|x|\to\infty}V_*(\Phi(x))=V_*(x_0)>0.
\]
This proves the existence of the asserted representative and completes the proof of Theorem~\ref{thm:nine}.

\end{proof}

	\section{Critical profile decomposition}\label{sec:critical-profile}
	
	\noindent
	\begin{revision}
The proof follows the critical strategy of \cite[Section~9]{FGL},
with the vector-valued fractional profile decomposition of
Lemma~\ref{lem:vector-profile}. We first extract profiles from the first
\(M\) eigenfunctions, where \(M\) counts the limiting Birman--Schwinger
eigenvalues at or above one. Rank allocation and the sharp mass bounds
then give \(M=N\). We treat profiles at the same center but different
scales by disjoint punctured neighborhoods, and use uniform convexity to
recover the potentials strongly.
\end{revision}
	
	\begin{proof}[Proof of Theorem~\ref{thm:seven}]
		\medskip
		\noindent\emph{Step 1: select the Birman--Schwinger eigenfunctions.}
		The sequence \(V_n\) is bounded in \(L^{q_c}\), so the compact operators \(\cT_{V_n}\) are uniformly bounded.
		\begin{revision}
The spectral bound needed here follows directly from the CLR inequality
\cite{FrankCLR}. If \(\mu_k(V)>0\), Proposition~\ref{prop:birman-schwinger}
gives, for every \(\varepsilon>0\),
\[
 k\leq N_-\!\left(( -\Delta)^s-
        \frac{1+\varepsilon}{\mu_k(V)}V\right)
 \leq L_{0,d,s}\frac{(1+\varepsilon)^{q_c}}{\mu_k(V)^{q_c}}
        \int V^{q_c}.
\]
Letting \(\varepsilon\searrow0\), and noting that the estimate is
trivial when \(\mu_k(V)=0\), yields the bound below.
A diagonal extraction therefore gives
\end{revision}
		\begin{equation*}
			\mu_k(V_n)\rightarrow\mu_k,
			\qquad
			\mu_k(V_n)
			\leq
			Ck^{-1/q_c}\|V_n\|_{q_c},
		\end{equation*}
		for every fixed \(k\).
		Since \(\mu_N(V_n)=1\), \(\mu_N=1\).
		The critical CLR bound therefore implies that there is a largest finite \(M\geq N\) with \(\mu_M\geq1\).
		Since \(\mu_k(V_n)\leq\mu_N(V_n)=1\) for \(k\geq N\), one has \(\mu_M=1\) and \(\mu_{M+1}<1\).
		Choose \(\dot H^s\)-orthonormal eigenfunctions \(f_{k,n}\), \(k=1,\ldots,M\), with
		\begin{equation*}
			\cT_{V_n}f_{k,n}
			=
			\mu_k(V_n)f_{k,n},
			\qquad
			(-\Delta)^sf_{k,n}
			=
			\mu_k(V_n)^{-1}V_nf_{k,n}.
		\end{equation*}
		
		\medskip
		\noindent\emph{Step 2: profile extraction and rank allocation.}
		Apply Lemma~\ref{lem:vector-profile} to \(F_n=(f_{1,n},\ldots,f_{M,n})\).
\begin{revision}
At least one profile is nonzero. Otherwise, the remainder conclusion
would imply \(\|f_{N,n}\|_{2_s^*}\to0\), contradicting the exact identity
\[
 1=\mu_N(V_n)=\int V_n|f_{N,n}|^2
 \leq\|V_n\|_{q_c}\|f_{N,n}\|_{2_s^*}^2.
\]
\end{revision}
		For each nonzero profile define
		\begin{equation*}
			f_{k,n}^{(j)}(x):=
			(t_n^{(j)})^{-(d-2s)/2}
			f_{k,n}\left(x_n^{(j)}+\frac{x}{t_n^{(j)}}\right),
		\end{equation*}
		so that \(f_{k,n}^{(j)}\weak f_k^{(j)}\) in \(\dot H^s\), and define the rescaled potentials
		\begin{equation*}
			V_n^{(j)}(x)
			:=
			(t_n^{(j)})^{-2s}
			V_n\left(x_n^{(j)}+\frac{x}{t_n^{(j)}}\right).
		\end{equation*}
		\begin{revision}
After a diagonal extraction over the at most countably many profiles,
\[
 V_n^{(j)}\weak V^{(j)}\quad\text{in }L^{q_c},
 \qquad V^{(j)}\geq0.
\]
The nonnegativity follows from weak closedness of the positive cone.
\end{revision}
		For \(\varphi\in C_c^\infty\), the rescaled eigenfunctions converge strongly in \(L^{q_c'}(\supp\varphi)\), since \(q_c'=d/(d-2s)<2_s^*\).
		Hence
		\begin{equation*}
			\int V_n^{(j)}f_{k,n}^{(j)}\varphi
			\rightarrow
			\int V^{(j)}f_k^{(j)}\varphi,
		\end{equation*}
		by splitting off the strong eigenfunction difference and then using weak \(L^{q_c}\)-convergence of the potentials.
		Passing to the limit in the weak eigenvalue equation gives
		\begin{equation*}
			(-\Delta)^sf_k^{(j)}
			=
			\mu_k^{-1}V^{(j)}f_k^{(j)}
		\end{equation*}
		for every nonzero profile component, initially with tests in
\(C_c^\infty(\R^d)\). H\"older's inequality and Sobolev embedding give
\[
 \left|\int V^{(j)}f_k^{(j)}\overline\varphi\right|
 \leq C_{d,s}\|V^{(j)}\|_{q_c}
       \|f_k^{(j)}\|_{\dot H^s}\|\varphi\|_{\dot H^s}.
\]
Thus both sides of the weak equation are continuous in the test function
with respect to the \(\dot H^s\)-norm. By density of \(C_c^\infty\) in
\(\dot H^s\), the equation holds for every \(\varphi\in\dot H^s\).
In particular, it may be tested with \(\varphi=f_k^{(j)}\), and it is
equivalent to \(\cT_{V^{(j)}}f_k^{(j)}=\mu_k f_k^{(j)}\).
		
		Writing \(L\) for the truncation index, the remainder estimate gives
		\begin{equation*}
			\lim_{L\to\infty}\limsup_{n\to\infty}
			\int V_n|r_{k,n}^{(L)}|^2
			\leq
			\sup_n\|V_n\|_{q_c}
			\lim_{L\to\infty}\limsup_{n\to\infty}
			\|r_{k,n}^{(L)}\|_{2_s^*}^2=0.
		\end{equation*}
		Every nonzero vector profile has at least one nonzero component.
		Testing its limiting equation against that component shows that \(V^{(j)}\neq0\) and that the corresponding \(\mu_k\geq1\) is an eigenvalue of \(\cT_{V^{(j)}}\).
		
		Let
		\begin{equation*}
			\mathcal E_j
			:=
			\one_{[1,\infty)}\bigl(\cT_{V^{(j)}}\bigr)\dot H^s,
			\qquad
			N_j:=\dim\mathcal E_j.
		\end{equation*}
		To justify the first rank inequality, fix a finite set \(\mathcal J\) of profiles and choose a \(\dot H^s\)-orthonormal eigenbasis \(g_1^{(j)},\ldots,g_{N_j}^{(j)}\) of \(\mathcal E_j\), with corresponding eigenvalues \(\nu_\ell^{(j)}\geq1\).
		Transport these functions to the original coordinates by
		\begin{equation*}
			g_{\ell,n}^{(j)}(x)
			:=
			(t_n^{(j)})^{(d-2s)/2}
			g_\ell^{(j)}\bigl(t_n^{(j)}(x-x_n^{(j)})\bigr).
		\end{equation*}
		Index the resulting finite family by \(\alpha=(j,\ell)\), and denote its Gram and potential-form matrices by
		\begin{equation*}
			(\mathsf G_n)_{\alpha\beta}
			:=\ipHs{g_{\ell,n}^{(j)}}{g_{r,n}^{(k)}},
			\qquad
			(\mathsf B_n)_{\alpha\beta}
			:=\int_{\R^d}V_n g_{\ell,n}^{(j)}
			\overline{g_{r,n}^{(k)}}.
		\end{equation*}
		Parameter orthogonality gives \(\mathsf G_n\to I\).
		On a diagonal block, a change of variables and the weak convergence \(V_n^{(j)}\weak V^{(j)}\) give
		\begin{equation*}
			(\mathsf B_n)_{(j,\ell),(j,r)}
			\rightarrow
			\int_{\R^d}V^{(j)}g_\ell^{(j)}\overline{g_r^{(j)}}
			=\nu_r^{(j)}\delta_{\ell r}.
		\end{equation*}
		If \(j\neq k\), H\"older's inequality gives
		\begin{equation*}
			\left|
			\int_{\R^d}V_n g_{\ell,n}^{(j)}
			\overline{g_{r,n}^{(k)}}
			\right|
			\leq
			\|V_n\|_{q_c}
			\|g_{\ell,n}^{(j)}g_{r,n}^{(k)}\|_{q_c'}
			\rightarrow0.
		\end{equation*}
		Here is a direct verification of the product convergence. Write
\(r=q_c'\), \(P=2r\), and \(a=(d-2s)/2=d/P\).
For fixed \(f,g\in L^P\), approximate them in \(L^P\) by bounded
compactly supported functions \(f_0,g_0\). H\"older and invariance of
the \(L^P\)-norm under \(T_n^{(j)}\) show that
\[
\begin{split}
 &\|(T_n^{(j)}f)(T_n^{(k)}g)
       -(T_n^{(j)}f_0)(T_n^{(k)}g_0)\|_r\\
 &\qquad\leq\|f-f_0\|_P\|g\|_P
       +\|f_0\|_P\|g-g_0\|_P,
\end{split}
\]
uniformly in \(n\). For the approximants, if the scale ratios stay bounded
above and below, parameter orthogonality forces their supports to be
disjoint for large \(n\). If, after interchanging the indices,
\(\lambda_n=t_n^{(k)}/t_n^{(j)}\to\infty\), change to the \(j\)-th frame.
Since \(2ar=d\), the product norm is unchanged by this common change
of frame and its \(r\)-th power is
\[
 \int |f_0(y)|^r\lambda_n^{ar}
       |g_0(\lambda_n(y-z_n))|^r\,dy
 \leq \|f_0\|_\infty^r\|g_0\|_r^r\lambda_n^{ar-d}
 =C\lambda_n^{-d/2}\to0,
\]
where \(z_n=t_n^{(j)}(x_n^{(k)}-x_n^{(j)})\).
Every subsequence has a further subsequence in one of these cases.
Letting the approximation errors tend to zero proves the asserted
convergence for the original fixed functions.
		Consequently
		\begin{equation*}
			\mathsf G_n^{-1/2}\mathsf B_n\mathsf G_n^{-1/2}
			\rightarrow
			\bigoplus_{j\in\mathcal J}
			\diag(\nu_1^{(j)},\ldots,\nu_{N_j}^{(j)}).
		\end{equation*}
		If \(r_{\mathcal J}:=\sum_{j\in\mathcal J}N_j\), the min--max principle therefore yields \(\mu_{r_{\mathcal J}}(V_n)\geq1-o(1)\).
		Since \(\mu_{M+1}(V_n)\to\mu_{M+1}<1\), it follows that \(r_{\mathcal J}\leq M\).
		Taking the supremum over finite \(\mathcal J\) gives
		\begin{equation*}
			\sum_jN_j\leq M.
		\end{equation*}
		Since every nonzero profile has \(N_j\geq1\), there are at most \(M\) profiles. Denote their total number by \(J\). The remainder is unchanged for \(L\geq J\), so
		\(\|r_{k,n}^{(J)}\|_{2_s^*}\to0\) as \(n\to\infty\).
		
		We still have to remove a possible homogeneous Sobolev remainder.
		Test the eigenvalue equation for \(f_{k,n}\) against \(r_{k,n}^{(J)}\).
		Profile orthogonality makes its \(\dot H^s\)-pairing with every extracted profile \(o(1)\), while H\"older and Sobolev give
		\begin{equation*}
			\left|\int V_nf_{k,n}\overline{r_{k,n}^{(J)}}\right|
			\leq
			\|V_n\|_{q_c}\|f_{k,n}\|_{2_s^*}
			\|r_{k,n}^{(J)}\|_{2_s^*}=o(1).
		\end{equation*}
		Since \(\mu_k(V_n)\to\mu_k\geq1\), the tested equation gives \(\|r_{k,n}^{(J)}\|_{\dot H^s}\to0\) as \(n\to\infty\), with \(J\) fixed.
		Applying \eqref{eq:vector-linear-combination-decoupling} to arbitrary coefficient vectors in \(\C^M\), and then using polarization, yields the Gram identity
		\begin{equation*}
			\sum_j\ipHs{f_k^{(j)}}{f_\ell^{(j)}}=\delta_{k\ell}.
		\end{equation*}
		For fixed \(j\), all vectors \(f_k^{(j)}\) belong to \(\mathcal E_j\).
		Define
		\begin{equation*}
			\mathsf A:\C^M\rightarrow\bigoplus_j\mathcal E_j,
			\qquad
			\mathsf A(c_1,\ldots,c_M)
			:=
			\left(\sum_{k=1}^M c_k f_k^{(j)}\right)_j.
		\end{equation*}
		The Gram identity gives, for every \(c\in\C^M\),
		\begin{equation*}
			\|\mathsf Ac\|_{\oplus_j\dot H^s}^2
			=
			\sum_{k,\ell=1}^M
			c_k\overline{c_\ell}
			\sum_j\ipHs{f_k^{(j)}}{f_\ell^{(j)}}
			=|c|^2.
		\end{equation*}
		Thus \(\mathsf A\) is an isometric injection, and hence \(M\leq\sum_jN_j\).
		Combining the two rank inequalities yields
		\begin{equation*}
			\sum_{j=1}^J N_j=M.
		\end{equation*}
		
		\medskip
		\noindent\emph{Step 3: punctured profile neighborhoods and the sharp chain.}
		We first construct the disjoint sets needed for the mass lower bound.
		This requires more than choosing large balls, because two profiles may have the same physical center and asymptotically different scales.
		
		After a further subsequence, every ratio \(t_n^{(k)}/t_n^{(j)}\)
has a limit in \([0,\infty]\). Call \(k\) narrower than \(j\) when this
ratio tends to infinity. For every such pair, extract again so that
\[
 z_{jk,n}:=t_n^{(j)}(x_n^{(k)}-x_n^{(j)})
\]
either converges to a finite point \(z_{jk}\), or has modulus tending
to infinity. These extractions are simultaneous because there are
only finitely many profiles. Let \(\mathcal P_j\) be the finite set
of finite limits \(z_{jk}\) associated with narrower profiles.

Choose \(\varepsilon_\ell\searrow0\) and radii \(R_\ell\nearrow\infty\)
such that, for every \(j\),
\[
 \int_{|x|>R_\ell}(V^{(j)})^{q_c}\leq\varepsilon_\ell.
\]
Absolute continuity of the finitely many integrals permits a decreasing
sequence \(\delta_\ell\searrow0\) such that
\[
 \int_{\bigcup_{z\in\mathcal P_j}B_{\delta_\ell}(z)}
       (V^{(j)})^{q_c}\leq\varepsilon_\ell
 \quad\hbox{for every }j.
\]
Define the compact sets
\[
 K_\ell^{(j)}:=\overline B_{R_\ell}(0)
       \setminus\bigcup_{z\in\mathcal P_j}B_{\delta_\ell}(z),
 \qquad
 \Omega_{n,\ell}^{(j)}:=x_n^{(j)}+(t_n^{(j)})^{-1}K_\ell^{(j)}.
\]
Then \(\int_{\R^d\setminus K_\ell^{(j)}}(V^{(j)})^{q_c}
\leq2\varepsilon_\ell\), and \(\one_{K_\ell^{(j)}}\to1\) pointwise
outside the finite set \(\mathcal P_j\), hence almost everywhere.

Fix \(\ell\) before taking \(n\to\infty\). For comparable scales,
parameter orthogonality gives
\(t_n^{(j)}|x_n^{(k)}-x_n^{(j)}|\to\infty\), whereas both radii
in the \(j\)-th frame remain bounded; their physical images are disjoint.
If \(k\) is narrower than \(j\), its image in the \(j\)-th frame lies in
\[
 z_{jk,n}+\frac{t_n^{(j)}}{t_n^{(k)}}\overline B_{R_\ell}(0).
\]
When \(|z_{jk,n}|\to\infty\), this lies outside
\(\overline B_{R_\ell}(0)\) for large \(n\). Otherwise, for large \(n\),
\[
 |z_{jk,n}-z_{jk}|+
 \frac{t_n^{(j)}}{t_n^{(k)}}R_\ell<\delta_\ell,
\]
so it lies entirely inside the hole deleted from \(K_\ell^{(j)}\).
Taking the maximum of the finitely many thresholds gives an index
\(n_0(\ell)\) beyond which all the \(\Omega_{n,\ell}^{(j)}\) are
pairwise disjoint. Thus the construction uses only the ordered limits
\(n\to\infty\) at fixed \(\ell\), followed by \(\ell\to\infty\).

		Scaling invariance and weak lower semicontinuity in the \(j\)-th frame give
		\begin{equation*}
			\liminf_{n\to\infty}
			\int_{\Omega_{n,\ell}^{(j)}}V_n^{q_c}
			\geq
			\int_{K_\ell^{(j)}}(V^{(j)})^{q_c}.
		\end{equation*}
		Summing over the eventually disjoint physical sets, then letting \(\ell\to\infty\), proves the required mass inequality
		\begin{equation}\label{eq:critical-profile-mass-lower-bound}
			\liminf_{n\to\infty}\int_{\R^d}V_n^{q_c}
			\geq
			\sum_{j=1}^J\int_{\R^d}(V^{(j)})^{q_c}.
		\end{equation}
		\begin{revision}
Combining \eqref{eq:critical-profile-mass-lower-bound} with the
definition of \(a_k\), the rank identity \(\sum_jN_j=M\), and
Lemma~\ref{lem:critical-constants} yields
\[
 a_N=\lim_{n\to\infty}\int V_n^{q_c}
 \geq\sum_j\int(V^{(j)})^{q_c}
 \geq\sum_j a_{N_j}\geq a_M\geq a_N.
\]
Thus all inequalities are equalities.
\end{revision}
		Hence every \(V^{(j)}\) minimizes the mass under the constraint \(\mu_{N_j}\geq1\).
		Its \(N_j\)-th Birman--Schwinger eigenvalue must equal one: if \(\mu_{N_j}(V^{(j)})>1\), replacing \(V^{(j)}\) by \(\mu_{N_j}(V^{(j)})^{-1}V^{(j)}\) would preserve admissibility and strictly decrease its \(L^{q_c}\)-mass.
		Therefore every \(V^{(j)}\) is an optimizer and
		\begin{equation*}
			\mu_{N_j}(V^{(j)})=1.
		\end{equation*}
		
		\medskip
		\noindent\emph{Step 4: exclusion of \(M>N\).}
		Assume \(M>N\).
		Recall \(a_k=k/\ell_{0,d,s}^{(k)}\).
		Equality in the preceding chain says
		\begin{equation*}
			a_N=\sum_{j=1}^Ja_{N_j}=a_M.
		\end{equation*}
		If some \(N_j=1\), Lemma~\ref{lem:critical-constants} and positivity of \(a_1\) give
		\begin{equation*}
			a_{M-1}\leq\sum_{\ell\neq j}a_{N_\ell}
			=a_N-a_1<a_N,
		\end{equation*}
		contrary to monotonicity \(a_{M-1}\geq a_N\).
		Hence every \(N_j\geq2\).
		Fix \(j_0\).
		Since \(\mu_{N_{j_0}-1}(V^{(j_0)})\geq1\), the variational definition of \(a_{N_{j_0}-1}\), Lemma~\ref{lem:critical-constants}, and \(M-1\geq N\) give
		\begin{align*}
			a_N
			&\leq a_{M-1}
			\leq a_{N_{j_0}-1}+\sum_{j\neq j_0}a_{N_j}\leq
			\frac{\displaystyle\int(V^{(j_0)})^{q_c}}
			{\mu_{N_{j_0}-1}(V^{(j_0)})^{q_c}}
			+\sum_{j\neq j_0}\int(V^{(j)})^{q_c}
			\leq a_N.
		\end{align*}
		All inequalities are equalities.
		Thus \(V^{(j_0)}\) is an optimizer for rank \(N_{j_0}-1\) and \(\mu_{N_{j_0}-1}(V^{(j_0)})=1\).
		Theorem~\ref{thm:nine} then gives \(\mu_{N_{j_0}}(V^{(j_0)})<1\), contradicting \(\mu_{N_{j_0}}(V^{(j_0)})=1\).
		Therefore
		\begin{equation*}
			M=N,
			\qquad
			\sum_{j=1}^J N_j=N.
		\end{equation*}
		
		\medskip
		\noindent\emph{Step 5: strong approximation by the profiles.}
		Equality of all critical masses gives
		\begin{equation*}
			\lim_{n\to\infty}\int V_n^{q_c}
			=
			\sum_{j=1}^J\int(V^{(j)})^{q_c}.
		\end{equation*}
		We recover the potentials while keeping track of the two limiting operations.
		Set
		\begin{equation*}
			U_n^{(j)}(x)
			:=
			(t_n^{(j)})^{2s}
			V^{(j)}\bigl(t_n^{(j)}(x-x_n^{(j)})\bigr).
		\end{equation*}
		Use the punctured compact sets \(K_\ell^{(j)}\) and their pairwise disjoint physical images \(\Omega_{n,\ell}^{(j)}\) constructed in Step~3.
		Their mass loss is \(O_J(\varepsilon_\ell)\), which is sufficient for the diagonal argument below.
		
		For fixed \(\ell\), weak convergence in the \(j\)-th frame and weak lower semicontinuity give
		\begin{equation*}
			\liminf_{n\to\infty}
			\int_{\Omega_{n,\ell}^{(j)}}V_n^{q_c}
			\geq
			\int_{K_\ell^{(j)}}(V^{(j)})^{q_c}.
		\end{equation*}
		Write \(b_j:=\int(V^{(j)})^{q_c}\) and \(b_{j,\ell}:=\int_{K_\ell^{(j)}}(V^{(j)})^{q_c}\).
		Since the neighborhoods are disjoint and \(\int V_n^{q_c}\to\sum_jb_j\), the preceding lower bounds also imply
		\begin{align*}
			\limsup_{n\to\infty}
			\int_{\R^d\setminus\bigcup_j\Omega_{n,\ell}^{(j)}}V_n^{q_c}
			&\leq\sum_j(b_j-b_{j,\ell}),\quad
			\limsup_{n\to\infty}
			\int_{\Omega_{n,\ell}^{(j)}}V_n^{q_c}
			\leq b_j+\sum_{k\neq j}(b_k-b_{k,\ell}).
		\end{align*}
		Thus the exterior mass is \(O(J\varepsilon_\ell)\), and the mass in the \(j\)-th neighborhood differs from \(b_{j,\ell}\) by at most \(O(J\varepsilon_\ell)\), after \(n\) is sufficiently large.
		
		Pass to a diagonal subsequence \(n=n_\ell\), and relabel it, so that these estimates hold and, for every \(j\),
		\begin{equation*}
			h_{j,\ell}
			:=
			\one_{K_\ell^{(j)}}V_{n_\ell}^{(j)}
			\weak V^{(j)}
			\quad\text{in }L^{q_c}.
		\end{equation*}
		For clarity, choose a countable dense family \((\psi_m)\) in
\(L^{q_c'}\). For each fixed \(\ell\), weak convergence in the profile
frames permits \(n_\ell\geq n_0(\ell)\), increasing in \(\ell\), such that
\[
 \left|\int\one_{K_\ell^{(j)}}
       (V_{n_\ell}^{(j)}-V^{(j)})\psi_m\right|\leq\ell^{-1}
 \quad(j=1,\ldots,J,\ m\leq\ell),
\]
and the preceding mass bounds hold with an additional error at most
\(\ell^{-1}\). Dominated convergence gives
\(\|(1-\one_{K_\ell^{(j)}})\psi_m\|_{q_c'}\to0\).
Uniform \(L^{q_c}\)-boundedness and density extend the resulting
convergence from \((\psi_m)\) to every dual test function, proving
the displayed weak convergence.
		The local mass bounds and \(b_{j,\ell}\to b_j\) give
		\begin{equation*}
			\|h_{j,\ell}\|_{q_c}^{q_c}\rightarrow b_j
			=\|V^{(j)}\|_{q_c}^{q_c}.
		\end{equation*}
		Since \(1<q_c<\infty\), uniform convexity of \(L^{q_c}\) now yields \(h_{j,\ell}\to V^{(j)}\) strongly.
		Since also \(\one_{K_\ell^{(j)}}V^{(j)}\to V^{(j)}\) strongly, returning to the physical variables gives
		\begin{equation*}
			\left\|
			\one_{\Omega_{n_\ell,\ell}^{(j)}}
			\bigl(V_{n_\ell}-U_{n_\ell}^{(j)}\bigr)
			\right\|_{q_c}
			\rightarrow0.
		\end{equation*}
		Moreover, the \(L^{q_c}\)-norm of \(V_{n_\ell}\) outside the union of the neighborhoods tends to zero, and
		\begin{equation*}
			\left\|
			\one_{\R^d\setminus\Omega_{n_\ell,\ell}^{(j)}}
			U_{n_\ell}^{(j)}
			\right\|_{q_c}^{q_c}
			=b_j-b_{j,\ell}\rightarrow0.
		\end{equation*}
		Decomposing the difference over the disjoint neighborhoods and their complement, and then using the triangle inequality, gives
		\begin{equation*}
			\left\|
			V_n-
			\sum_{j=1}^J
			(t_n^{(j)})^{2s}
			V^{(j)}\bigl(t_n^{(j)}(\cdot-x_n^{(j)})\bigr)
			\right\|_{q_c}
			\to0.
		\end{equation*}
		The equality of constants is exactly the displayed formula in the theorem.
	\end{proof}
	
	\begin{proof}[Proof of Corollary~\ref{cor:critical-compactness}]
		Theorem~\ref{thm:seven} gives a partition \(N=N_1+\cdots+N_J\) and the identity \(a_N=\sum_{j=1}^Ja_{N_j}\).
		If \(J\geq2\), set \(K=N_1\).
		Repeated use of \eqref{eq:critical-subadditivity} gives \(a_{N-K}\leq\sum_{j=2}^Ja_{N_j}\); hence \eqref{eq:critical-strict-subadditivity} would imply
		\begin{equation*}
			a_N<a_K+a_{N-K}
			\leq\sum_{j=1}^Ja_{N_j}=a_N,
		\end{equation*}
		a contradiction.
		Thus \(J=1\), and the strong \(L^{q_c}\) approximation in Theorem~\ref{thm:seven} gives the asserted compactness modulo translations and dilations.
		If \(\ell_{0,d,s}^{(N)}>\max_{1\leq k<N}\ell_{0,d,s}^{(k)}\), then for every \(1\leq K<N\),
		\begin{equation*}
			\frac{K}{\ell_{0,d,s}^{(K)}}
			+\frac{N-K}{\ell_{0,d,s}^{(N-K)}}
			>
			\frac{K}{\ell_{0,d,s}^{(N)}}
			+\frac{N-K}{\ell_{0,d,s}^{(N)}}
			=\frac{N}{\ell_{0,d,s}^{(N)}},
		\end{equation*}
		which is \eqref{eq:critical-strict-subadditivity}.
		This is the argument of \cite[Corollary~8]{FGL} after the critical profile decomposition has been established.
	\end{proof}
	
	\begin{revision}
\begin{proof}[Proof of Theorem~\ref{thm:rank-two}]
Write \(q=q_c\), \(P=2_s^*\), and \(S=S_{d,s}\).
Lemma~\ref{lem:critical-constants} gives \(a_2\leq2a_1=2S^q\).
We prove the reverse inequality and nonattainment simultaneously.

Let \(0\leq V\in L^q\) with \(\mu:=\mu_2(V)>0\).
Choose a real eigenfunction \(h\) for \(\mu_1(V)\).
The inequality \(\||h|\|_{\dot H^s}\leq\|h\|_{\dot H^s}\) shows,
by the Rayleigh principle, that \(|h|\) is also a maximizer of the
first Rayleigh quotient. We may therefore take \(h\geq0\) with
\(\mathcal T_Vh=\mu_1(V)h\).
H\"older, HLS and uniqueness in \(\dot H^s\), as in Step~1 of
Proposition~\ref{prop:critical-decay}, give
\[
 h=\mu_1(V)^{-1}I_{2s}(Vh)>0\quad\text{almost everywhere}.
\]
Here \(Vh\not\equiv0\), since
\(\int Vh^2=\mu_1(V)\|h\|_{\dot H^s}^2>0\).
Choose a nonzero real eigenfunction \(u\) for \(\mu\) with
\(\langle u,h\rangle_{\dot H^s}=0\). Self-adjointness gives
\[
 \int Vuh=0,\qquad \int Vu^2=\mu\|u\|_{\dot H^s}^2>0.
\]
Consequently \(u\) has both signs on \(\{V>0\}\), so
\(u^+,u^-\neq0\).

Test \(( -\Delta)^su=\mu^{-1}Vu\) with \(u^+\) and \(-u^-\).
The strict cross-term identity \eqref{eq:nodal-cross-term} yields
\[
 \mu^{-1}\int V(u^\pm)^2
 =\|u^\pm\|_{\dot H^s}^2
    -\langle u^+,u^-\rangle_{\dot H^s}
 >\|u^\pm\|_{\dot H^s}^2.
\]
On the other hand,
\[
 \int V(u^\pm)^2
 \leq\left(\int_{\{\pm u>0\}}V^q\right)^{1/q}\|u^\pm\|_P^2
 \leq S^{-1}\left(\int_{\{\pm u>0\}}V^q\right)^{1/q}
             \|u^\pm\|_{\dot H^s}^2.
\]
Dividing by the nonzero Sobolev norms gives
\[
 \int_{\{\pm u>0\}}V^q>(\mu S)^q,
 \qquad
 \int V^q>2S^q\mu_2(V)^q.
\]
Thus every potential with \(\mu_2(V)\geq1\) has mass strictly
larger than \(2S^q\). Taking the infimum and combining with
subadditivity gives \(a_2=2a_1\), while the strict inequality excludes
attainment. Amplitude normalization gives the equivalent nonattainment
of the quotient defining \(\ell_{0,d,s}^{(2)}\).

For a normalized optimizing sequence, Theorem~\ref{thm:seven} gives
\(\sum_jN_j=2\). A single profile would be a rank-two optimizer,
which has just been excluded. Hence \(J=2\) and \(N_1=N_2=1\).
The rank-one classification in Section~1.1 identifies both profiles as
translates and dilations of \(Q_{d,s}\). Absorbing these fixed
transformations into their parameters preserves parameter orthogonality
and gives the asserted approximation.
Finally, strong \(L^q\)-convergence of any translated and dilated
subsequence would preserve the mass and \(\mu_2\), since
\[
 |\mu_2(V)-\mu_2(W)|\leq\|\mathcal T_V-\mathcal T_W\|
 \leq C_{d,s}\|V-W\|_q.
\]
Its limit would be a rank-two optimizer, a contradiction.
\end{proof}

\begin{proof}[Proof of Corollary~\ref{cor:strict-mass-gap}]
Set \(a_0=0\). We first show that
\[
 a_{m-1}<a_m\quad\text{whenever }a_m\text{ is attained}.
\]
For \(m=1\), this is the positivity of \(a_1\).
For \(m\geq2\), let \(W\) be a normalized optimizer for \(a_m\).
If \(a_{m-1}=a_m\), then \(W\) also minimizes the rank-\((m-1)\)
mass, because \(\mu_{m-1}(W)\geq\mu_m(W)=1\).
Amplitude scaling forces \(\mu_{m-1}(W)=1\); otherwise the mass
could be reduced. Theorem~\ref{thm:nine}, applied at rank \(m-1\),
would give \(\mu_m(W)<1\), a contradiction.

Now apply Theorem~\ref{thm:seven} to a normalized minimizing sequence
for \(a_{N+1}\). It provides attained masses \(a_{N_j}\) with
\[
 \sum_{j=1}^JN_j=N+1,\qquad \sum_{j=1}^Ja_{N_j}=a_{N+1}.
\]
Fix \(j_0\). Subadditivity, omitting a zero summand if \(N_{j_0}=1\),
and the preceding strict inequality give
\[
 a_N\leq a_{N_{j_0}-1}+\sum_{j\ne j_0}a_{N_j}
 <\sum_{j=1}^Ja_{N_j}=a_{N+1}.
\]
This proves strict increase without assuming attainment at rank \(N+1\).

For every nonnegative \(V\in L^{q_c}\), amplitude normalization in
\eqref{eq:critical-minimization} gives
\[
 \mu_{N+1}(V)^{q_c}\leq\frac{1}{a_{N+1}}\int V^{q_c};
\]
the inequality is immediate when \(\mu_{N+1}(V)=0\).
For \(\int V_n^{q_c}\to a_N\), taking the upper limit yields
\[
 \limsup_n\mu_{N+1}(V_n)
 \leq\left(\frac{a_N}{a_{N+1}}\right)^{1/q_c}<1.
\]
Applying the same inequality to a normalized optimizer proves the
last assertion.
\end{proof}
\end{revision}


\begin{thebibliography}{99}
		\phantomsection
		\addcontentsline{toc}{section}{References}
		
		\bibitem{CY26}
L. Chen and L. Yang,
\emph{Optimizers for the finite-rank fractional Lieb--Thirring inequality},
arXiv:2610.21716v1 (2026),
\url{https://arxiv.org/abs/2610.21716v1}.

\bibitem{Cwi77}
		M.~Cwikel,
		\emph{Weak type estimates for singular values and the number of bound states of
			Schr\"odinger operators},
		Ann. of Math. (2) 106 (1977), 93--100.
		
		\bibitem{FrankCLR}
		R.~L. Frank,
		\emph{Cwikel's theorem and the CLR inequality},
		J. Spectr. Theory 4 (2014), 1--21.
		
		\bibitem{FrankSurvey}
		R.~L. Frank,
		\emph{The Lieb--Thirring inequalities: recent results and open problems},
		in \emph{Nine Mathematical Challenges---An Elucidation},
		Proc. Sympos. Pure Math. 104, American Mathematical Society, 2021, 45--86.
		
		\bibitem{FGL21a}
		R.~L. Frank, D.~Gontier, and M.~Lewin,
		\emph{The nonlinear Schr\"odinger equation for orthonormal functions II:
			application to Lieb--Thirring inequalities},
		Comm. Math. Phys. 384 (2021), 1783--1828.
		
		\bibitem{FGL}
		R.~L. Frank, D.~Gontier, and M.~Lewin,
		\emph{Optimizers for the finite-rank Lieb--Thirring inequality},
		Amer. J. Math. 147 (2025), 503--560.
		
		\bibitem{FLW}
		R.~L. Frank, A.~Laptev, and T.~Weidl,
		\emph{Schr\"odinger Operators: Eigenvalues and Lieb--Thirring Inequalities},
		Cambridge Studies in Advanced Mathematics 200, Cambridge University Press, 2023.
		
		\bibitem{Gerard}
		P.~G\'erard,
		\emph{Description du d\'efaut de compacit\'e de l'injection de Sobolev},
		ESAIM Control Optim. Calc. Var. 3 (1998), 213--233.
		
		\bibitem{Kato}
		T.~Kato,
		\emph{Perturbation Theory for Linear Operators},
		Classics in Mathematics, Springer-Verlag, 1995.
		
		\bibitem{Lie76}
		E.~H. Lieb,
		\emph{Bounds on the eigenvalues of the Laplace and Schr\"odinger operators},
		Bull. Amer. Math. Soc. 82 (1976), 751--753.
		
		\bibitem{Lieb83}
E.~H. Lieb,
\emph{Sharp constants in the Hardy--Littlewood--Sobolev and related inequalities},
Ann. of Math. (2) 118 (1983), 349--374.

\bibitem{LiebSeiringer}
		E.~H. Lieb and R.~Seiringer,
		\emph{The Stability of Matter in Quantum Mechanics},
		Cambridge University Press, 2010.
		
		\bibitem{LT75}
		E.~H. Lieb and W.~Thirring,
		\emph{Bound for the kinetic energy of fermions which proves the stability of matter},
		Phys. Rev. Lett. 35 (1975), 687--689; erratum, 1116.
		
		\bibitem{LT76}
		E.~H. Lieb and W.~Thirring,
		\emph{Inequalities for the moments of the eigenvalues of the Schr\"odinger
			Hamiltonian and their relation to Sobolev inequalities},
		in \emph{Studies in Mathematical Physics}, Princeton University Press, 1976,
		269--303.
		
		\bibitem{PalatucciPisante}
		G.~Palatucci and A.~Pisante,
		\emph{Improved Sobolev embeddings, profile decomposition, and
			concentration--compactness for fractional Sobolev spaces},
		Calc. Var. Partial Differential Equations 50 (2014), 799--829.
		
		\bibitem{RS75}
		M.~Reed and B.~Simon,
		\emph{Methods of Modern Mathematical Physics. II. Fourier Analysis,
			Self-Adjointness},
		Academic Press, 1975.
		
		\bibitem{RS78}
		M.~Reed and B.~Simon,
		\emph{Methods of Modern Mathematical Physics. IV. Analysis of Operators},
		Academic Press, 1978.
		
		\bibitem{Roz72}
		G.~V. Rozenblum,
		\emph{Distribution of the discrete spectrum of singular differential operators},
		Dokl. Akad. Nauk SSSR 202 (1972), 1012--1015.
		
		\end{thebibliography}
	
\end{document}